# DENSITIES, SPECTRAL DENSITIES AND MODALITY[1]

By P. Laurie Davies and Arne Kovac

*Universität Duisburg-Essen*

This paper considers the problem of specifying a simple approximating density function for a given data set $(x_1, \ldots, x_n)$. Simplicity is measured by the number of modes but several different definitions of approximation are introduced. The taut string method is used to control the numbers of modes and to produce candidate approximating densities. Refinements are introduced that improve the local adaptivity of the procedures and the method is extended to spectral densities.

**1. Contents.** In Section 1.1 we formulate the density problem in terms of obtaining the simplest density which is an adequate approximation for the given data. The taut string method of Davies and Kovac (2001) is adapted to the density problem and is used for producing candidate densities of increasing complexity. The difficulties of the density problem are discussed in Section 2. Section 3 contains a more detailed account of the application of the taut string method to the density problem. The asymptotics of the procedure on appropriate test beds are discussed in Section 4. A refinement based on cell occupancy frequencies which increases local sensitivity is described in Section 5. Section 5.4 compares the taut string method with kernel estimators in a small simulation study. Finally, Section 6 describes the application of the taut string methodology to the problem of spectral densities.

1.1. *The density problem.* Given a sample $\mathbf{x}_n = (x_1, \ldots, x_n)$ of size $n$, we consider the problem of specifying a distribution $F$ with the smallest number of modes such that the resulting model of i.i.d. random variables

Received April 2001; revised June 2003.

[1]Supported in part by Sonderforschungsbereich 475, University of Dortmund and DFG Grants 237/2-1 and 237/2-2 on time series analysis.

$\mathbf{X}_n^F = (X_1^F, \ldots, X_n^F)$ with common distribution $F$ is an adequate approximation for the data $\mathbf{x}_n$.

We use different concepts of approximation, one of which is the following. Let $E_n$,

$$E_n(x) = \frac{1}{n} \sum_{i=1}^{n} \{x_i \leq x\},$$

denote the empirical distribution of the data $\mathbf{x}_n$ and $F_n$ the empirical distribution function of $n$ i.i.d. random variables $\mathbf{X}_n^F$ with common distribution $F$. The Kolmogorov metric $d_{\text{KO}}$ is defined by

$$d_{\text{KO}}(F, G) = \sup\{x : |F(x) - G(x)|\}.$$

The i.i.d. model with distribution $F$ will be regarded as an adequate approximation to the data $\mathbf{x}_n$ if

(1.1) $$d_{\text{KO}}(E_n, F) \leq \text{qu}(n, \alpha, d_{\text{KO}}),$$

where $\text{qu}(n, \alpha, d_{\text{KO}})$ denotes the $\alpha$-quantile of the random variable $d_{\text{KO}}(F_n, F)$ which is independent of $F$ for continuous $F$. This gives rise to the Kolmogorov problem.

PROBLEM 1.1 (Kolmogorov problem). Determine the smallest integer $k_n$ for which there exists a density $f^n$ with $k_n$ modes and whose distribution $F^n$ satisfies

(1.2) $$d_{\text{KO}}(E_n, F^n) \leq \text{qu}(n, \alpha, d_{\text{KO}}).$$

We note that the problem is well posed: for any data set $\mathbf{x}_n$ it has a solution. We have posed the problem in terms of approximation so that no assumptions regarding the "true" data generating mechanism are required or made.

The Problem 1.1 is formulated in terms of the smallest number of modes required for an adequate approximation. A detailed theoretical discussion of such one-sided problems is given by Donoho (1988); one of his examples is that of modality of nonparametric densities and spectral densities. His paper also raises interesting questions about statistical inference involving objects whose very existence cannot be shown, an example being the "underlying density" for the data. We avoid such problems by phrasing the paper in terms of approximation.

Hartigan and Hartigan (1985) and Hartigan (2000) construct tests for the modality of a density function. They are based on the Kolmogorov distance of the nearest mixture of uniform distributions to the data and are discussed in more detail below.

Hengartner and Stark (1995) also make use of the Kolmogorov ball to determine nonparametric confidence bounds for densities subject to an upper bound for the number of modes. In the particular case of monotone or unimodal densities the width of their bounds on appropriate test beds is $(\log n/n)^{1/3}$, which agrees with the results given in this paper. It seems that their bounds become difficult to calculate for more than one mode as the complexity is given as $\binom{n}{l}$ where $l$ is the number of local extremes. The main differences with respect to the work of Hengartner and Stark (1995) are as follows:

(i) we provide an explicit density but no bounds,
(ii) neither the number of modes nor even an upper bound is specified in advance,
(iii) the algorithmic complexity of our method is $O(n)$ independently of the number of modes.

1.2. *The taut string methodology.* The basic methodology we use for producing densities is the taut string methodology. Taut strings were first used in the context of monotonic regression: the greatest convex minorant of the integrated data is a taut string and its derivative is precisely the monotone increasing least squares approximation. This is described in Barlow, Bartholomew, Bremner and Brunk (1972), who were the first to use the phrase "taut string." We refer also to Leurgans (1982). The first use of the taut string which goes beyond the monotone case and which explicitly deals with modality is in Hartigan and Hartigan (1985), where it is referred to as the "stretched string." Hartigan and Hartigan (1985) introduced their DIP test for unimodality which is based on the closest (in the Kolmogorov metric) unimodal distribution to the empirical distribution function of the data. Based on the work of Hartigan and Hartigan (1985), Davies (1995) used the taut string method to produce candidate densities of low modality to approximate data. Mammen and van de Geer (1997) employed the taut string in the nonparametric regression problem. They considered a penalized least squares problem where the penalty is the total variation of the approximating function. The solution is the basic taut string confined to a tube centered at the integrated data. Mammen and van de Geer (1997) gave a detailed description of the taut string but did not mention the connection with modality. Hartigan (2000) recently proposed a generalization of the DIP test to an arbitrary number of modes. It is based on the Kolmogorov distance between the empirical distribution and the nearest distribution consisting of a mixture of uniform distributions with at most $m$ modes. This is calculated using a taut string. Hartigan examines for each antimode of a taut string approximation the supremum distance between the empirical distribution function and a monotone density on a "shoulder interval" including the antimode. Finally, Davies and Kovac (2001) used the taut string methodology to control the number of local extremes of a nonparametric approximation to a data set. They also introduced the idea of local squeezing and residual driven tube widths, which greatly increases the precision and flexibility of the taut string methodology.

1.3. *Smoothness.* The taut string methodology produces densities which are piecewise constant and therefore not even continuous. Smoothness will not be a consideration in this paper but we point out that techniques for smoothing such functions have been developed. The idea is to obtain the smoothest density subject to shape and deviation constraints taken from the taut string. We refer to Metzner (1997), Löwendick and Davies (1998) and Majidi (2003).

1.4. *Previous work.* Much work has been done on the problem of density estimation. One of the most popular methods is that of kernel smoothing. We refer to Watson (1964), Nadaraya (1965), Silverman (1986), Sheather and Jones (1991), Wand and Jones (1995), Sain and Scott (1996) and Simonoff (1996) and the references given therein. The main problem here is the determination of appropriate global or local bandwidths. A further approach is based on wavelets. We refer to Donoho, Johnstone, Kerkyacharian and Picard (1996), Herrick, Nason and Silverman (2000) and Vidakovic (1999), Chapter 7. Mixtures of densities have been considered in the Bayesian framework by Richardson and Green (1997) and Roeder and Wasserman (1997). Other Bayesian methods are to be found in Verdinelli and Wasserman (1998).

None of the above approaches is directly concerned with modality. For example, the non-Bayesian theory is generally based on integrated squared error or some similar loss function. In spite of this, methods are often judged by their ability to identify peaks in the data as in Loader (1999) and Herrick, Nason and Silverman (2000). Work directly concerned with modality has been done by Müller and Sawitzki (1991) using their concept of excess mass. Their ideas have been extended to multidimensional distributions by Polonik (1995a, b, 1999). Hengartner and Stark (1995) used the Kolmogorov ball centered at the empirical distribution function to obtain nonparametric confidence bounds for shape restricted densities. Another way of controlling modality is that of mode testing. We refer to Good and Gaskins (1980), Silverman (1986), Hartigan and Hartigan (1985) and Fisher, Mammen and Marron (1994).

**2. The difficulties of the density problem.** Obtaining adequate approximate densities is a special case of nonparametric regression. Whereas nonparametric regression is usually concerned with the size of the dependent variable, the density problem is concerned with measuring the degree of closeness of the design points. In spite of a formal similarity, this is the more difficult problem and it may explain the modesty evident in the literature on densities. The difficulties may be illustrated by three data sets each of a sample size of $n = 500$. The first was generated using the standard normal distribution, the second using the uniform distribution on $[0, 1]$ and the third using the so-called claw distribution which is the following mixture of five normal distributions:

$$0.5\mathcal{N}(0,1) + 0.1\sum_{i=0}^{4}\mathcal{N}(i/2 - 1, 0.1).$$

This density will also be referred to as $N5$ (see Section 3.1). It is one of ten introduced by Marron and Wand (1992) to study the performance of different density methods. For each data set we calculated a kernel estimate with a global bandwidth which was chosen to be as small as possible subject to the estimate having the same modality as the density. Similarly for the taut string method we took the Kolmogorov ball to be as small as possible subject to the estimate having the same modality as the density. The results are shown in Figure 1.

The kernel method performs very well on the sample from the normal distribution but the approximation to the uniform density is poor. It can only be improved by using a smaller bandwidth which then introduces superfluous modes. The approximation to the claw density is even worse. Only three peaks are correctly identified; the remaining two peaks are in the tails near $-2$ and $3$, where the claw density does not have a peak. An explanation of this behavior can be found in Hartigan (2000), who discusses the relationship between the peaks and bandwidth for kernel estimates.

The taut string method produces excellent approximations in all three cases. In particular, all five peaks of the claw density are correctly identified. The open problem is to produce an automatic procedure for the taut string method which will give good approximations on these and other test beds without knowledge of the number of modes. In the case of nonparametric regression such an automatic procedure is available and is reminiscent of hard thresholding for wavelets [Davies and Kovac (2001)]. Unfortunately, there seems to be no equivalent for densities and it is this which makes the density problem so difficult.

## 3. Taut strings, Kuiper metrics and densities.

3.1. *Test densities.* As part of the evaluation of the procedures to be defined below, we consider their performance on test beds defined by distributions. For the sake of convenient reference we list here the distributions we consider. $\mathcal{N}(\mu, \sigma^2)$ refers to the normal distribution with mean $\mu$ and variance $\sigma^2$.

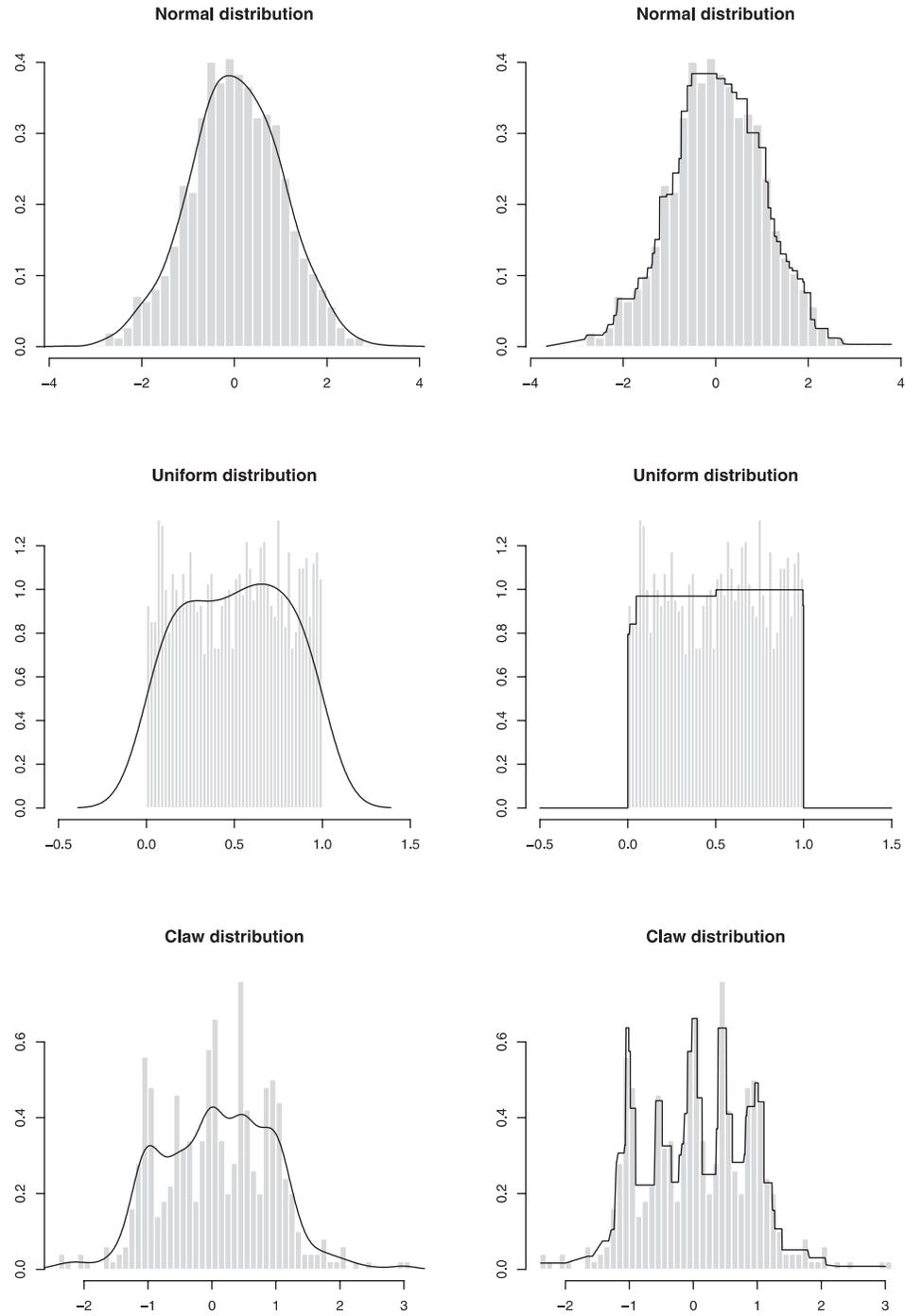

FIG. 1. *Normal, uniform and claw density. The panels show kernel and taut string approximations using the smallest bandwidth that retains the correct modality.*

| | |
|---|---|
| $U$ | the uniform distribution on $[0,1]$ |
| $N1$ | the standard normal distribution |
| $S$ | the slash distribution, defined as $\mathcal{N}(0,1)/\mathcal{U}(0,1)$ [see Morgenthaler and Tukey (1991)] |
| $N2$ | the mixture $0.5\mathcal{N}(0,1)+0.5\mathcal{N}(3,1)$ |
| $N4$ | the mixture $0.8\mathcal{N}(0,3)+0.015\mathcal{N}(8,0.02)+0.015\mathcal{N}(9,0.02)+0.17\mathcal{N}(15,0.2)$ |
| $N5$ | the claw distribution $0.5\mathcal{N}(0,1)+0.1\sum_{i=0}^{4}\mathcal{N}(i/2-1,0.1)$ |
| $N10\_5$ | the mixture $0.1\sum_{i=1}^{10}\mathcal{N}(5i-5,1)$ |
| $N10\_10$ | the mixture $0.1\sum_{i=1}^{10}\mathcal{N}(10i-5,1)$ |

3.2. *Taut strings.* We give a short description of the taut string method. A thorough analysis of properties of the taut string can be found in Hartigan (2000). Further details and an algorithm of complexity $O(n)$ are given by Davies and Kovac (2001).

Consider a sample $\mathbf{x}_n$ and form the ordered sample $\mathbf{x}_{(n)}=(x_{(1)},\ldots,x_{(n)})$. For a given $\varepsilon>0$ we consider the Kolmogorov tube $T(E_n,\varepsilon)$ centered at the empirical distribution $E_n$ and of radius $\varepsilon>0$

$$T(E_n,\varepsilon)=\left\{G:G\text{ monotone }\sup_x|G(x)-E_n(x)|\leq\varepsilon\right\}.$$

Imagine now a taut string taking the value of 0 at $x_{(1)}$ and 1 at $x_{(n)}$ and constrained to lie within the Kolmogorov tube. Such a string is shown in the right panels of Figure 2 for two different values of $\varepsilon$. The taut string defines a function $S_n$ on the interval $[x_{(1)},x_{(n)}]$. Although $S_n$ depends on $E_n$ and $\varepsilon$, we suppress this dependency to relieve the burden on the notation. We denote the density of $S_n$ by $s_n$. It is defined as the left-hand side derivative of $S_n$ except at the smallest data point $x_{(1)}$ where we use the right-hand side derivative. The left panels of Figure 2 show histograms of the data with the corresponding densities $s_n$ superimposed.

The taut string is a spline with knots at the points at which it touches the lower or upper boundaries of the Kolmogorov tube. The taut string has the following properties [see Davies and Kovac (2001) and Mammen and van de Geer (1997)]:

 (i) $S_n$ is monotonic increasing and linear between knots.
 (ii) $s_n$ is nonnegative and piecewise constant between knots.
(iii) $s_n$ has the minimum modality of all functions whose integral over $[x_{(1)},x_{(n)}]$ lies in $T(E_n,\varepsilon)$ and satisfies the end point conditions.
(iv) $S_n$ switches from the upper boundary $E_n+\varepsilon$ to the lower boundary $E_n-\varepsilon$ at points where $s_n$ has a local maximum.
 (v) $S_n$ switches from the lower boundary $E_n-\varepsilon$ to the upper boundary $E_n+\varepsilon$ at points where $s_n$ has a local minimum.

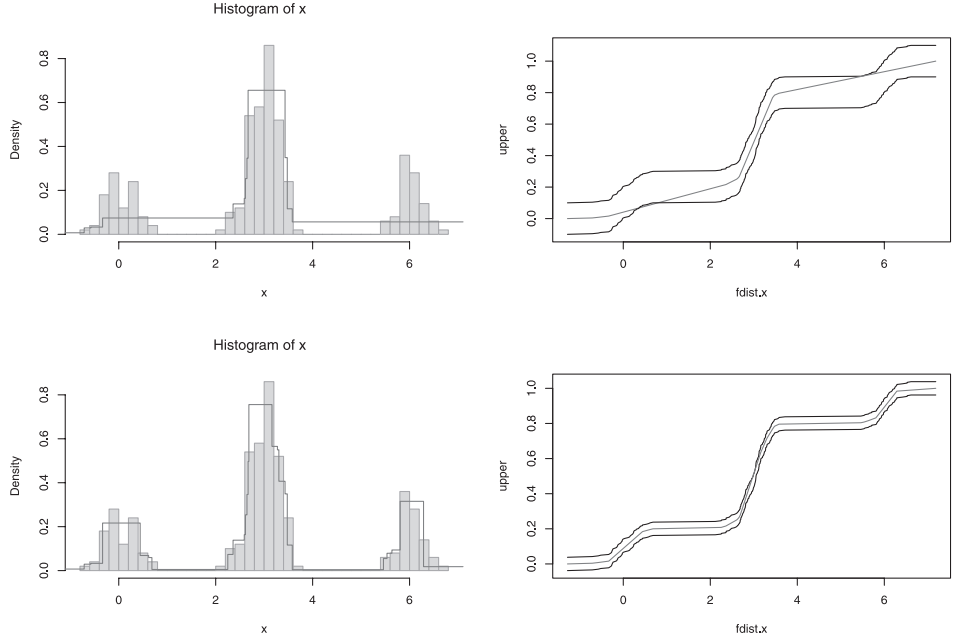

Fig. 2. *These figures illustrate the taut string method applied to a sample of mixture of normal distributions with two different tubewidths. The right column shows the tubes and the taut strings while the left column shows histograms of the data and the corresponding densities of the taut string.*

(vi) If $\xi_j$ and $\xi_{j+1}$ are consecutive knots on the same boundary, then on the interval $(\xi_j, \xi_{j+1}]$

$$s_n(x) = \frac{|\{i : \xi_j < x_i \leq \xi_{j+1}\}|}{n(\xi_{j+1} - \xi_j)}. \tag{3.1}$$

It is property (iii) which is of importance and allows control of the number of modes. If consecutive knots $\xi_j$ and $\xi_{j+1}$ are on opposite boundaries, then it follows from (iv) and (v) above that (3.1) must be replaced by

$$s_n(x) = \frac{|\{i : \xi_j < x_i \leq \xi_{j+1}\}| \pm 2\varepsilon}{n(\xi_{j+1} - \xi_j)} \tag{3.2}$$

with a minus sign at local maxima and a plus sign at local minima. This means that the derivative underestimates local maxima and overestimates local minima. In an earlier version of this paper we followed Davies and Kovac (2001) and modified string $\tilde{S}_n$ by setting

$$\tilde{S}_n(\xi_j) = E_n(\xi_j) \qquad \text{at all knots } \xi_j \tag{3.3}$$

and linear in between. The corresponding derivative $\tilde{s}_n$ satisfies

$$(3.4) \quad \tilde{s}_n(\xi_j) = \frac{|\{i : \xi_j < x_i \leq \xi_{j+1}\}|}{n(\xi_{j+1} - \xi_j)} \qquad \text{between the knots } \xi_j \text{ and } \xi_{j+1}.$$

This modification has no effect on the modality and in general produces more pronounced peaks. More by good luck than by good thinking, the authors fortunately noticed that much improved results can be obtained by *not* modifying the taut string in this manner. The reason is that this alteration causes both the taut string and the empirical distribution to have the same mass on intervals defining local extremes. Below we shall use Kuiper metrics which are defined by those intervals where the difference is greatest. The idea is that differences in distributions with different peaks should be greatest on intervals defining peaks. Modifying the taut string as in (3.4) nullifies this effect. Nevertheless, the *final* density, which is returned by the procedure, is modified in this manner.

3.3. *Data analysis.* Even without an automatic procedure, the taut string can be used as a data analytical tool. If the radius of the Kolmogorov tube is monotonically decreased, then the number of modes of the derivative of the taut string increases monotonically. It is therefore possible to specify the number of modes of the approximate density. Figure 3 shows this for the same sample as used for Figure 1. The densities of Figure 3 can be interpreted as histograms with an automatic choice of the number of bins and the bin widths. To measure the performance of the taut string procedure, we simulated samples of different sizes from the claw distribution and squeezed the tube as far as possible consistent with the density having five peaks. A peak is classified as being correctly identified if the midpoint of the interval defining a peak differs by less than 0.15 from the position of the nearest peak of the claw density. Figure 4 shows the number of correctly identified peaks as a function of sample size.

It shows that the taut string method is extremely good at finding peaks. For samples of size 200, the five peaks will be correctly identified in over 80% of the cases. This in a sense confirms Loader (1999), who, on the basis of theoretical results of Marron and Wand (1992), claims that for samples of size $n = 193$ the claws should be detectable. The problem we now address is the difficult one of defining an automatic procedure with a similar performance.

3.4. *An automatic procedure.* The following theorem is an immediate consequence of the properties of the taut string listed above.

THEOREM 3.1. *The derivative $s_n$ of the taut string constrained to lie in the tube $T(E_n, \mathrm{qu}(n, \alpha, d_{\mathrm{KO}}))$ is a solution of the Kolmogorov density problem.*

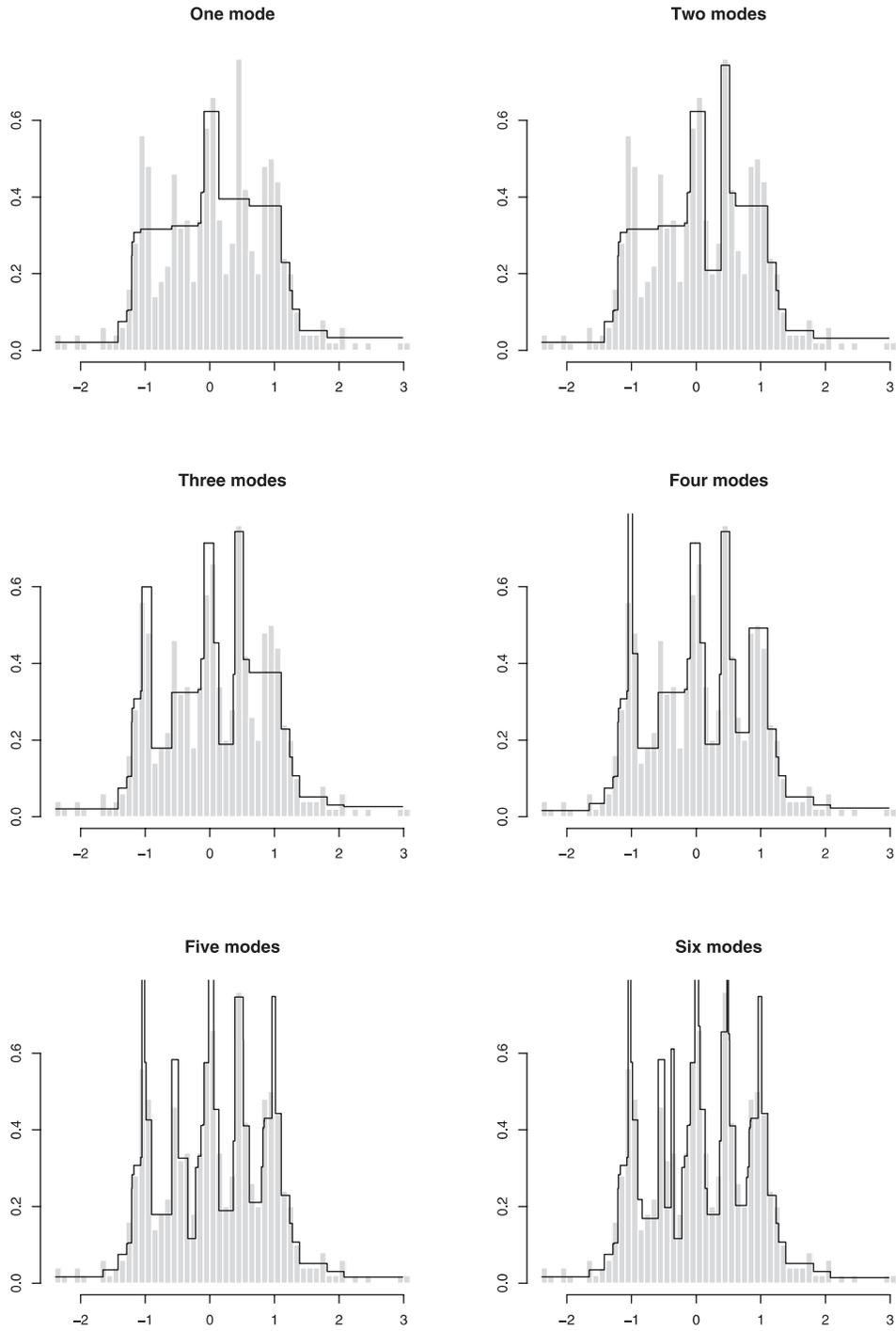

Fig. 3. *Six taut string estimates of a sample of the claw distribution with increasing number of modes.*

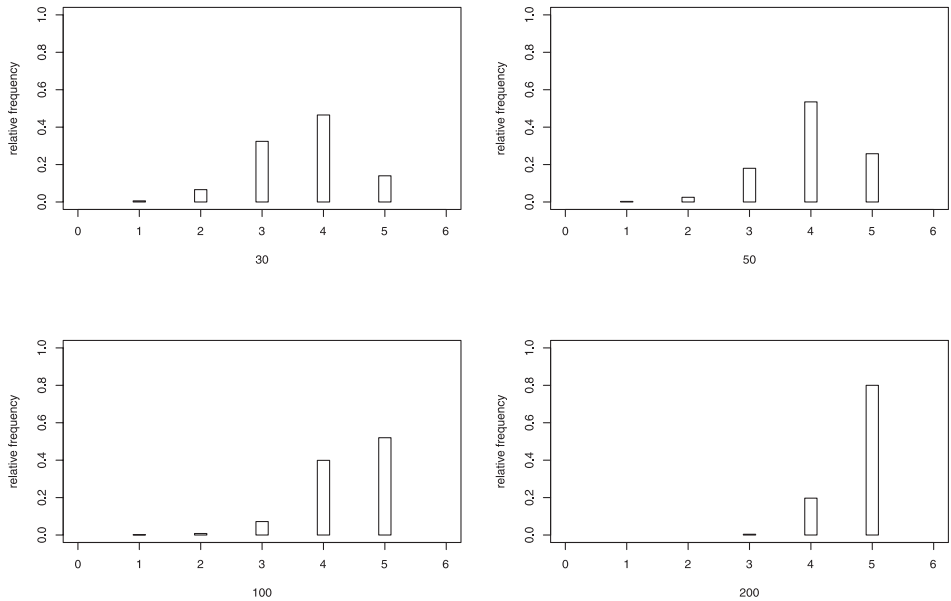

Fig. 4. *Five-modal taut string*: *number of correctly identified peaks of the claw density as a function of sample size.*

For finite $n$ the values of $\mathrm{qu}(n,\alpha,d_{\mathrm{KO}})$ can be obtained by simulation. In the limit $\sqrt{n}\,\mathrm{qu}(n,\alpha,d_{\mathrm{KO}})$ tends to the corresponding quantile of

$$\max_{0\leq t\leq 1} B_0(t) - \min_{0\leq t\leq 1} B_0(t),$$

where $B_0$ denotes a Brownian bridge and for which an explicit expression exists [Dudley (1989)].

The solution of the Kolmogorov density problem therefore defines an automatic procedure based on the taut string and its performance can be evaluated on different test beds. If we do this on an i.i.d. test bed, that is, with data of the form $X_1(F),\ldots,X_n(F)$ where $F$ has a $k$-modal density function $f$, then it is clear that the taut string density $s_n$ will have at most $k$ modes with probability at least $\alpha$. This follows on noting that $F$ lies in the tube with probability $\alpha$ and that in this case $s_n$ has at most as many modes as $f$. In particular, if $k=1$, we have the following.

THEOREM 3.2. *Let $X_1(F),\ldots,X_n(F)$ be an i.i.d. sample with common unimodal distribution $F$ and let $s_n$ be the solution of the Kolmogorov density problem* (1.2). *Then*

(3.5) $$\mathbb{P}(s_n\ unimodal) \geq \alpha.$$

A simulation was performed to investigate the performance of the procedure with $\alpha = 0.9$ and the corresponding tube width $1.245/\sqrt{n}$ on test beds defined by the distributions listed in Section 3.1. The results are shown in Table 1. It is clear that for a unimodal distribution the modality is correctly estimated with probability at least 0.9 in accordance with Theorem 3.2. Indeed the actual probability greatly exceeds 0.9 as all simulations resulted in exactly one peak. The results for the other distributions are, in contrast, disappointing. Asymptotically the modality will be correctly estimated with probability at least 0.9 but the rate of convergence is clearly very slow. We now try and obtain an improved procedure in two ways. First we note that the choice of $\text{qu}(n, \alpha, d_{\text{KO}})$ for the radius of the tube means that a probability of at least $\alpha$ is guaranteed for all unimodal test beds. If we provisionally accept that the uniform distribution is a poor model for most data sets, then we may accept a worse performance for the uniform distribution in return for enhanced performances for other distributions. Silverman (1986) and Müller and Sawitzki (1991) argue in a similar vein. The second way of gaining an improved performance is to use a generalized Kuiper metric rather than the Kolmogorov metric. Kuiper metrics consider the differences in probability over a fixed number of disjoint intervals and are therefore better at detecting modality.

3.5. *Calibrating unimodality.* To implement the first way of improving performance, let $\text{qu}(n, \alpha, F, 1, d_{\text{KO}})$ denote the $\alpha$-quantile of the Kolomogorov distance of the closest unimodal distribution (given by the taut string) to the empirical distribution $F_n$ of $n$ i.i.d. random variables with common distribution $F$. We have the following theorem.

THEOREM 3.3. *Let $X_1(F), \ldots, X_n(F)$ be an i.i.d. sample with common unimodal distribution $F$ and empirical distribution $F_n$. Let $s_n$ be the deriva-*

TABLE 1
*The procedure using the* 0.9*-quantile of the Kolmogorov metric. The numbers give the percentage of simulations in which the correct modality was obtained. The numbers in parentheses give the mean absolute deviation from the correct modality. The results are based on* 1000 *simulations*

| Dist. | $U$ | $S$ | $N1$ | $N2$ | $N4$ | $N5$ | $N10\_5$ | $N10\_10$ |
|---|---|---|---|---|---|---|---|---|
| 100 | 100 (0) | 100 (0) | 100 (0) | 0 (1) | 0 (2.34) | 0 (4) | 0 (9) | 0 (9) |
| 500 | 100 (0) | 100 (0) | 100 (0) | 0 (1) | 0 (2) | 0 (4) | 0 (9) | 0 (8.6) |
| 1000 | 100 (0) | 100 (0) | 100 (0) | 0 (1) | 0 (2) | 0 (4) | 0 (9) | 0 (7.9) |
| 5000 | 100 (0) | 100 (0) | 100 (0) | 50 (0.5) | 0 (2) | 0 (4) | 0 (8.3) | 100 (0) |
| 10000 | 100 (0) | 100 (0) | 100 (0) | 100 (0) | 0 (2) | 66 (0.4) | 99 (0.01) | 100 (0) |

*tive of the string $S_n$ through the tube $T(F_n, \mathrm{qu}(n, \alpha, F, 1, d_{\mathrm{KO}}))$. Then*

(3.6) $$\mathbb{P}(s_n \text{ unimodal}) = \alpha.$$

Clearly

$$\mathrm{qu}(n, \alpha, F, 1, d_{\mathrm{KO}}) \leq \mathrm{qu}(n, \alpha, d_{\mathrm{KO}}),$$

but it is not clear whether

$$\sup_{F \text{ unimodal}} \mathrm{qu}(n, \alpha, F, 1, d_{\mathrm{KO}}) = \mathrm{qu}(n, \alpha, d_{\mathrm{KO}}).$$

We point out that the uniform distribution does not maximize $\mathrm{qu}(n, \alpha, F, 1, d_{\mathrm{KO}})$ [Hartigan and Hartigan (1985)]. We now take $F = U$ to be the uniform distribution on the basis that it is not an adequate approximation for most data sets and set $\alpha = 0.5$. This means that on uniform test beds the modality will be correctly determined with probability 0.5. The uniform distribution has the advantage that the asymptotics of the quantiles $\mathrm{qu}(n, \alpha, U, 1, d_{\mathrm{KO}})$ can be calculated. We have

(3.7) $$\lim_{n \to \infty} \sqrt{n}\, \mathrm{qu}(n, \alpha, U, 1, d_{\mathrm{KO}}) = \mathrm{qu}(\alpha, B_0),$$

where $\mathrm{qu}(\alpha, B_0)$ denotes the $\alpha$-quantile of the random variable

(3.8) $$\min_H \sup_x |B_0(x) - H(x)|,$$

where the function $H : [0,1] \to \mathbb{R}$ is convex on $[0, t_H]$ and concave on $[t_H, 1]$ for some $t_H, 0 \leq t_H \leq 1$. Simulations show that the 0.5-quantile of (3.8) is 0.432. A correction for finite $n$ gives

$$\mathrm{qu}(n, 0.5, U, 1, d_{\mathrm{KO}}) = 0.43/\sqrt{n} - 0.64/n,$$

with a percentage error (based on simulations) of at most 0.0045. Table 2 shows the results. We see that the performance for the Gaussian test bed is hardly impaired. On the claw test bed we note that the performance for $n = 1000$ is now comparable to that of the simple Kolmogorov quantile for $n = 10000$.

If we apply the same idea to the normal distribution, then heuristic arguments indicate that

$$\lim_{n \to \infty} \sqrt{n}\, \mathrm{qu}(n, \alpha, \mathcal{N}(0,1), 1, d_{\mathrm{KO}}) = 0$$

but we have no exact asymptotic result. The same argument goes through for any sufficiently smooth density. If true, this implies that if we use a cut-off point for the size of the Kolmogorov ball which is bounded below by some constant multiple of $1/\sqrt{n}$, then the modality will be consistently estimated. We do not pursue this idea any further.

TABLE 2
*The procedure based on the 0.5-quantile of the Kolmogorov distance of the closest unimodal distri- bution to a uniform sample. The numbers give to the nearest integer the percentage of simulations in which the correct modality was obtained. The numbers in parentheses give the mean absolute deviation from the correct modality correct to one decimal place. The results are based on* 1000 *simulations*

| Dist. | $S$ | $N1$ | $N2$ | $N4$ | $N5$ | $N10\_5$ | $N10\_10$ |
|---|---|---|---|---|---|---|---|
| 100   | 100 (0) | 100 (0) | 22 (0.8) | 0 (2)    | 0 (3.8)  | 0 (8)    | 0 (3.8) |
| 500   | 100 (0) | 100 (0) | 78 (0.2) | 0 (2)    | 1 (2.5)  | 0 (5.5)  | 1 (2.5) |
| 1000  | 100 (0) | 100 (0) | 95 (0)   | 0 (2)    | 43 (0.7) | 27 (1.1) | 43 (0.7) |
| 5000  | 100 (0) | 100 (0) | 99 (0)   | 48 (0.6) | 100 (0)  | 100 (0)  | 100 (0) |
| 10000 | 100 (0) | 100 (0) | 100 (0)  | 100 (0)  | 100 (0)  | 100 (0)  | 100 (0) |

3.6. *Kuiper metrics.* Suppose that the density $s_n$ of the taut string is unimodal. Part of the description of the taut string $S_n$ given in Section 3.2 is that it switches from the upper bound to the lower bound at each maximum. Consider now the Kuiper metric $d_{\mathrm{KU}}$ defined by

$$(3.9) \quad d_{\mathrm{KU}}(F,G) = \sup\{a < b : |(F(b) - F(a)) - (G(b) - G(a))|\}.$$

It follows from the above that if $d_{\mathrm{KO}}(E_n, S_n) = \varepsilon$ and $s_n$ is unimodal, then $d_{\mathrm{KU}}(E_n, S_n) = 2\varepsilon$. The $\alpha$-quantile $\mathrm{qu}(n, \alpha, d_{\mathrm{KU}})$ of $d_{\mathrm{KU}}(F_n, F)$ is independent of $F$ for continuous $F$ and is less than twice the $\alpha$-quantile of $d_{\mathrm{KO}}(F_n, F)$. This suggests that the Kuiper metric is more appropriate for unimodality than the Kolmogorov metric. To demonstrate this we firstly define the Kuiper problem.

PROBLEM 3.1 (Kuiper density problem). Determine the smallest integer $k_n$ for which there exists a density $f^n$ with $k_n$ modes and whose distribution $F^n$ satisfies

$$d_{\mathrm{KU}}(E_n, F^n) \leq \mathrm{qu}(n, \alpha, d_{\mathrm{KU}}).$$

Suppose now that $F^n$ is a unimodal distribution which solves the Kuiper density problem. Let

$$\varepsilon_1 = \max\{x : F^n(x) - E_n(x)\} \quad \text{and}$$
$$\varepsilon_2 = \max\{x : G(x) - F^n(x)\}.$$

As $d_{\mathrm{KU}}(E_n, F^n) = \varepsilon_1 + \varepsilon_2 = \mathrm{qu}(n, \alpha, d_{\mathrm{KU}})$, it follows by shifting $F^n$ by an amount $\frac{1}{2}|\varepsilon_2 - \varepsilon_1|$ that the solution of the Kolmogorov problem with $\varepsilon = \frac{1}{2}\mathrm{qu}(n, \alpha, d_{\mathrm{KU}})$ is also unimodal. As $\frac{1}{2}\mathrm{qu}(n, \alpha, d_{\mathrm{KU}}) < \mathrm{qu}(n, \alpha, d_{\mathrm{KO}})$, this implies that if the solution of the Kuiper density problem for a given $\alpha$ is unimodal, so is the solution of the Kolmogorov problem for the same $\alpha$.

To cover the case of multimodality we define the Kuiper metric $d_{\text{KU}}^\kappa$ of order $\kappa$ by

$$(3.10) \quad d_{\text{KU}}^\kappa(F,G) = \max\left\{\sum_1^\kappa |(F(b_j) - F(a_j)) - (G(b_j) - G(a_j))|\right\},$$

where the maximum is taken over all $a_j, b_j$ with

$$a_1 \le b_1 \le a_2 \le b_2 \le \cdots \le a_\kappa \le b_\kappa.$$

Again the distribution of $d_{\text{KU}}^\kappa(F_n, F)$ is independent of $F$ for continuous $F$. If we denote the $\alpha$-quantile by $\text{qu}(n, \alpha, d_{\text{KU}}^\kappa)$, we can formulate the $\kappa$-Kuiper problem.

PROBLEM 3.2 ($\kappa$-Kuiper density problem). Determine the smallest integer $k_n$ for which there exists a density $f^n$ with $k_n$ modes and whose distribution $F^n$ satisfies

$$d_{\text{KU}}^k(E_n, F^n) \le \text{qu}(n, \alpha, d_{\text{KU}}^\kappa).$$

If the density $s_n$ of the taut string has $k$ modes, then for the Kuiper metric $d_{\text{KU}}^{2k-1}$ of order $2k-1$ we have

$$d_{\text{KU}}^{2k-1}(E_m, S_n) = (2k-1)\varepsilon.$$

This follows on noting that the string switches boundaries at each of the $k$ local maxima of $s_n$ and also at the $k-1$ local minima. As

$$\text{qu}(n, \alpha, d_{\text{KU}}^{2k-1}) < (2k-1)\text{qu}(n, \alpha, d_{\text{KO}}),$$

this indicates that the Kuiper metric $d_{\text{KU}}^{2k-1}$ is more efficacious when the data exhibit $k$ modes. We have no simple algorithm for solving the $\kappa$-Kuiper problem so we use the strategy of Davies and Kovac (2001) and decrease the radius $\varepsilon$ of the Kolmogorov tube gradually until

$$d_{\text{KU}}^{2k-1}(E_n, S_n) \le \text{qu}(n, \alpha, d_{\text{KU}}^{2k-1}).$$

For large $n$ approximations to $\text{qu}(n, \alpha, d_{\text{KU}}^\kappa)$ are available using the weak convergence result

$$\sqrt{n} d_{\text{KU}}^\kappa(F_n, F) \Rightarrow \max\left\{\sum_1^\kappa |B_0(b_j) - B_0(a_j)|\right\},$$

where $B_0$ denotes the standard Brownian bridge on $[0,1]$ and

$$a_1 < b_1 < a_2 < b_2 < \cdots < a_\kappa < b_\kappa.$$

The distribution of $\max\{|B_0(b) - B_0(a)|\}$ corresponding to the unimodal case $k = 1$ is known [e.g., Dudley (1989), Proposition 12.3.4.]. Sufficiently accurate quantiles for finite $n$ and for the other asymptotic cases may be obtained by simulations. Best results are obtained if $\kappa$ is related to the modality $k$ of the test bed by $\kappa = 2k - 1$. In practice a default value of $\kappa$ is required and we use $\kappa = 19$.

We combine the $\kappa$-Kuiper metric with the ideas of Section 3.5. Let $\text{qu}(n, \alpha, F, 1, d_{\text{KU}}^\kappa)$ denote the $\alpha$-quantile of the $\kappa$-Kuiper distance of the closest unimodal distribution to the empirical distribution $F_n$ of $n$ i.i.d. random variables with common distribution $F$. We use the string $S_n$ as the closest unimodal distribution. If $F$ is the uniform distribution of $[0, 1]$, then we have again a $1/\sqrt{n}$ asymptotic. For example, for $\kappa = 19$ and $\alpha = 0.5$, simulations showed that

$$\text{qu}(n, 0.5, U, 1, d_{\text{KU}}^{19}) \approx 8.12/\sqrt{n} - 30.32/n^{1.04}$$

is a good approximation.

The results shown in Table 3 confirm the claim that the Kuiper metric with $\kappa = 2k - 1$ performs best on test beds with $k$ modes. Thus the procedure based on the $d_{\text{KU}}^3$ metric is best for the bimodal distribution $N2$, that based on the $d_{\text{KU}}^9$ metric is best for the five-modal claw density $N5$, while that based on the $d_{\text{KU}}^{19}$ metric is best for the two ten-modal distributions $N10\_5$ and $N10\_10$. None of the procedures performs well for the four-modal $N4$ distribution. This is because it has two very concentrated but lower power peaks situated at the points 8 and 9. For this distribution global squeezing of the Kolmogorov tube will only work for large sample sizes. In small samples when the tube is sufficiently narrow to pick up the lower power peaks, it will have already caused peaks to appear at other points. This is shown by Table 4. For the sample sizes shown the tube was squeezed to give just four peaks and it was then checked if the four peaks were the correct ones. Table 4 gives the percentage of cases when this was the case. Thus even for a sample of size 2000, the correct peaks were only found in 80% of the cases. The problem is related to that of detecting low power peaks in nonparametric regression. In Davies and Kovac (2001) the problem was solved using local squeezing. In Section 5 we introduce a form of local squeezing for densities which is based on cell occupancy frequencies.

TABLE 3
*Results for the procedures using the* 0.5-*quantile of the closest unimodal distribution in the Kuiper metrics based on* 3, 9 *and* 19 *intervals. The numbers give the percentage of simulations in which the correct modality was obtained. The numbers in parentheses give the mean absolute deviation from the correct modality. The results are based on* 1000 *simulations with sample sizes of* 250 *and* 500

| Dist.    | $S$     | $N1$    | $N2$      | $N4$    | $N5$      | $N10\_5$  | $N10\_10$ |
|----------|---------|---------|-----------|---------|-----------|-----------|-----------|
| $n=250$  |         |         |           |         |           |           |           |
| $k=3$    | 99 (0)  | 96 (0)  | 67 (0.3)  | 0 (2)   | 0 (2.9)   | 0 (6.7)   | 38 (0.8)  |
| $k=9$    | 100 (0) | 99 (0)  | 59 (0.4)  | 0 (1.9) | 20 (1.5)  | 0 (3.4)   | 95 (0)    |
| $k=19$   | 100 (0) | 96 (0)  | 53 (0.5)  | 1 (1.9) | 20 (1.5)  | 0 (1.0)   | 99 (0)    |
| $n=500$  |         |         |           |         |           |           |           |
| $k=3$    | 99 (0)  | 99 (0)  | 90 (0.1)  | 0 (2)   | 10 (1.7)  | 0 (3.9)   | 100 (0)   |
| $k=9$    | 100 (0) | 99 (0)  | 74 (0.3)  | 1 (1.9) | 70 (0.3)  | 50 (0.6)  | 100 (0)   |
| $k=19$   | 100 (0) | 99 (0)  | 66 (0.3)  | 2 (1.9) | 57 (0.5)  | 97 (0)    | 100 (0)   |

TABLE 4
*Results of global squeezing for the four-modal distribution N4. The Kolmogorov tube was squeezed to give exactly four peaks. The numbers give the percentage of simulations in which these were the correct peaks. The results are based on* 1000 *simulations*

| $n$ | 500 | 1000 | 2000 | 4000 |
|---|---|---|---|---|
|  | 3 | 23 | 81 | 99 |

3.7. *Discrete data.* So far we have looked for an approximation to the data in the form of a Lebesgue density. However, at little cost we can extend the methodology to integer-valued data which typically arise from counts. Suppose the data set $\mathbf{x}_n = (x_1, \ldots, x_n)$ contains only $N$ different values $t_1 < t_2 < \cdots < t_N$. We look for an approximation in terms of $N$ probabilities $p_j = \mathbb{P}(X = t_j)$, $j = 1, \ldots, N$, where the random variable $X$ has support $t_1 < t_2 < \cdots < t_N$. Let $e_1, \ldots, e_N$ be the empirical frequencies of the $t_j$ in the data and consider the cumulative sums

$$E_j = \sum_{i=1}^{j} e_i$$

and the tube constructed by linear interpolation of the points $(j/N, E_j)$, $j = 0, \ldots, N$. Differentiating yields an approximation of $p_1, \ldots, p_N$. This procedure corresponds to the taut string algorithm in the regression context [Davies and Kovac (2001)] with time points $t_1, \ldots, t_n$ and with observations $e_1, \ldots, e_n$. Our default procedure uses the $\kappa$-Kuiper metric with $\kappa = 9$ and $\alpha = 0.5$. We point out that this radius is conservative for discrete data, but we do not pursue this any further. Other forms of approximation can be accommodated without much difficulty. An example is shown in Figure 5 where the discrete taut string method was applied to 1200 observations from a mixture of three Poisson distributions,

$$0.25P(2) + 0.5P(8) + 0.25P(21).$$

The other situation is where repeated values occur not because of the na-ture of the data (counting) but because of rounding. The rounding of data is very common and it can cause difficulties when looking for an approximation based on Lebesgue densities. To see the difficulties assume that some data point $x$ is observed $k$ times. Depending on the exact implementation of the taut string algorithm, two problems may occur. If the tube is centered around the empirical distribution function and the tube width is smaller than $k/2n$, the derivative of the taut string at $x$ will be $\infty$. If, on the other hand, the tube is constructed by linear interpolation of the empirical distribution function, then the empirical mass at $x$ of $k/n$ is spread over the interval $[x_l, x]$, where $x_l$ is the largest data point smaller than $x$. To overcome these problems we propose the following. Let $\varepsilon$ be the precision or cut-off error which we set to $\varepsilon = 10^{-3}\mathrm{MAD}(\mathbf{x}_n)$, where MAD denotes the median absolute deviation. We construct a modified data set $\tilde{x}_1, \ldots, \tilde{x}_n$, where the identical observations at $x$ are equally spread over the interval $[x - \varepsilon/2, x + \varepsilon/2]$. To be precise, we replace $x_{(j+1)} = x_{(j+2)} = \cdots = x_{(j+k)}$ by

$$\tilde{x}_{j+i} = x + \varepsilon\left(-\frac{1}{2} + \frac{1}{2k} + \frac{i-1}{k}\right)$$

for $i = 1, \ldots, k$. The taut string method described above is then applied to $\tilde{x}$ instead of $x$.

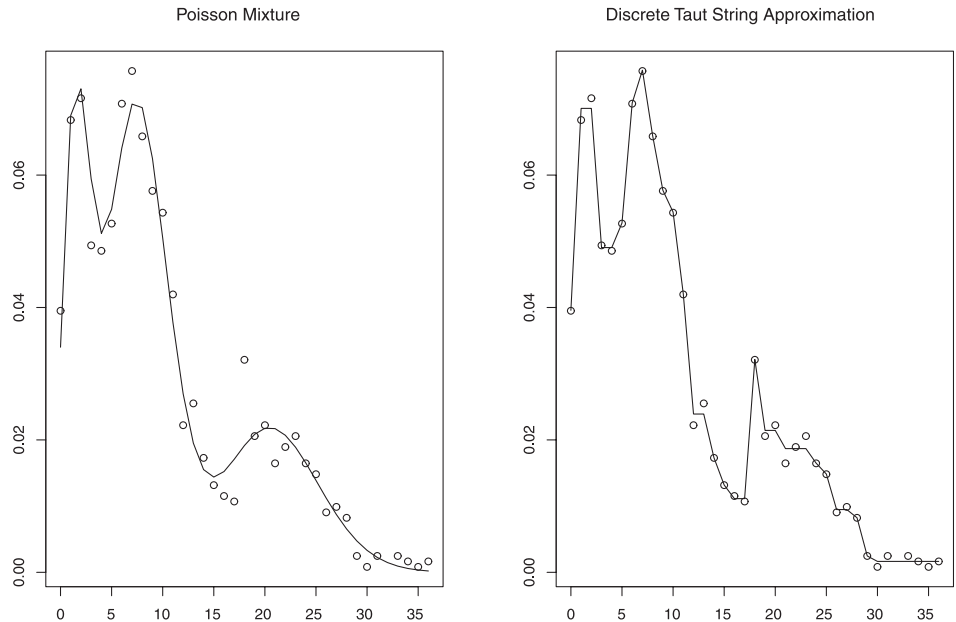

Fig. 5. *Discrete data. The left panel shows the density function of the mixture of three Poisson distributions and the frequencies of a sample of size* 1200. *The discrete taut string approximation is shown in the right panel.*

**4. Asymptotics on test beds.** The asymptotic behavior of the taut string may be analyzed on appropriate test beds. It turns out that asymptotically the modality is correctly estimated and that the optimal rate of convergence is attained except in small intervals containing the local extremes of the density $f$.

We denote the modality of the derivative of the taut string in the supremum tube $T(F_n, C/\sqrt{n})$ by $k_n^C$. The taut string based on the radius $C/\sqrt{n}$ will be denoted by $S_n^C$ with derivative $s_n^C$. We write $I_i^e(n, C)$, $1 \leq i \leq k_n^C$, for the intervals where $s_n^C$ attains its local extreme values and denote the midpoints of these intervals by $t_i^e(n, C)$, $1 \leq i \leq k_n^C$. The length of an interval $I$ will be denoted by $|I|$.

THEOREM 4.1. *Let $f$ be a $k$-modal density function on $\mathbb{R}$ such that*

$$\min_{g, (k-1)\text{-modal}} |F(x) - G(x)| > 0.$$

*Then we have, for all $\delta > 0$,*

$$\lim_{C \to \infty} \liminf_{n \to \infty} \mathbb{P}\bigg( \{k_n^C = k\} \cap \bigg\{ \max_{1 \leq i \leq k_n^C} |I_i^e(n, C)| \leq \delta \bigg\}$$

$$\cap \bigg\{ \max_{1 \leq i \leq k_n^C} |t_i^e(n, C) - t_j^e| \leq \delta \bigg\} \bigg) = 1.$$

In the following $A$ denotes a generic constant which depends only on $f$ and whose value may differ from appearance to appearance.

THEOREM 4.2. *Assume that*

(i) *$f$ has a compact support on $[0, 1]$,*
(ii) *$f$ has exactly $k$ local extreme values at the points $0 < t_1^e < \cdots < t_k^e < 1$,*
(iii) *$f$ has a bounded second derivative $f^{(2)}$ which is nonzero at the $k$ local extremes,*
(iv) *$f^{(1)}(t) = 0$ only for $t \in \{t_1^e, \ldots, t_k^e\}$.*

*Then the following statements hold:*

(a)

$$\lim_{C \to \infty} \liminf_{n \to \infty} \mathbb{P}(t_i^e \in I_i^e(n, C), 1 \leq i \leq k) = 1.$$

(b) *For every $C_1 < 6$ and $C_2 > 12$,*

$$\lim_{C \to \infty} \liminf_{n \to \infty} \mathbb{P}\bigg( |I_i^e(n, C)| \bigg( \frac{\sqrt{n} |f^{(1)}(t_i^e)|}{C} \bigg)^{1/3} \in [C_1^{1/3}, C_2^{1/3}], 1 \leq i \leq k \bigg) = 1.$$

(c) *Let $\xi_j^{n,C}$ be the knots of the taut string $S_n^C$ and denote*

$$m(n,C) = \max\left\{\xi_{j+1}^{n,C} - \xi_j^{n,C} : \xi_j^{n,C}, \xi_{j+1}^{n,C} \in (0,1) \setminus \bigcup_1^k I_i^e(n,C)\right\}.$$

*For some constant $A$ only depending on $f$, we have*

$$\lim_{C \to \infty} \liminf_{n \to \infty} \mathbb{P}\left(m(n,C) \leq \left(A|f^{(1)}(x_j)|^{-2/3}\left(\frac{\log n}{n}\right)^{1/3}\right)\right) = 1.$$

(d) *Denote*

$$M(n,C) = \left[A\left(\frac{\log n}{n}\right)^{1/3}, 1 - A\left(\frac{\log n}{n}\right)^{1/3}\right] \setminus \bigcup_i^n I_i^e(n,C).$$

*Then for some constant $A$ only depending on $f$, we have*

$$\lim_{C \to \infty} \liminf_{n \to \infty} \mathbb{P}\left(\max_{t \in M(n,C)} |f(t) - f_n^C(t)| \leq \left(A|f^{(1)}(t)|^{1/3}\left(\frac{\log n}{n}\right)^{1/3}\right)\right) = 1.$$

(e) *For some constants $A_1$ and $A_2$ only depending on $f$, we have*

$$\lim_{C \to \infty} \liminf_{n \to \infty} \mathbb{P}\left(\max_{t \in \bigcup_1^n I_i^e(n,C)} |f(t) - f_n^C(t)| \leq AC^{2/3} n^{-1/3}\right) = 1.$$

Part (d) of the theorem shows that, bounded away from the local extrema, the taut string density attains the optimal rate of convergence up to a logarithmic factor. The proofs follow the lines of Davies and Kovac (2001) and we omit them.

## 5. Cell occupancy frequencies and local squeezing.

5.1. *Cell occupancy frequencies.* The multiresolution procedure of Davies and Kovac (2001) is based on comparing the residuals of some regression function with those of Gaussian white noise. The comparison is based on the means on intervals which form a multiresolution scheme. A similar idea can be applied to the density problem. A distribution $F$ is an adequate model for the data $\mathbf{x}_n = (x_1, \ldots, x_n)$ if the transformed data

$$\mathbf{u}_n = F(\mathbf{x}_n) = (F(x_1), \ldots, F(x_n))$$

looks like an i.i.d. sample of size $n$ from the uniform distribution on $[0,1]$. This is done by comparing the frequencies

$$w_{jk}^n = |\{l : k2^{-j} < u_l \leq (k+1)2^{-j}\}|, \qquad 0 \leq k \leq 2^j, 1 \leq j \leq m,$$

with those to be expected from i.i.d. uniform random variables. The maximum resolution level $m$ is taken to be the smallest integer such that $n \leq 2^m$.

Suppose that $U_1, \ldots, U_n$ are independently and uniformly distributed on $[0, 1]$. Then
$$W_{jk}^n = |\{l : k2^{-j} < U_l \leq (k+1)2^{-j}\}|$$
is binomially distributed $b(n, 2^{-j})$. For given $\alpha$ we define the upper bounds $\upsilon_{j,k}^n(\alpha)$ by

(5.1) $$\upsilon_j^n(\alpha) = \min\left\{\upsilon : \mathbb{P}(Z_j^n \geq \upsilon) \leq \frac{1-\alpha}{2n}\right\},$$

where $Z_j^n$ satisfies the binomial distribution $b(n, 2^{-j})$. It follows that
$$\mathbb{P}(W_{jk}^n < \upsilon_j^n(\alpha), 1 \leq k \leq 2^j, 1 \leq j \leq n) \geq \alpha.$$
Lower bounds $\lambda_{j,k}^n(\alpha)$ can be derived similarly. This gives rise to the following problem.

PROBLEM 5.1 (Cell occupancy problem). Determine the smallest integer $k_n$ for which there exists a density $f^n$ with $k_n$ modes and whose distribution $F^n$ is such that the cell frequencies $w_{j,k}^n$ satisfy

(5.2) $$\lambda_j^n(\alpha) \leq w_{j,k}^n \leq \upsilon_j^n(\alpha),$$

where the $\upsilon_{j,k}^n(\alpha)$ are given by (5.1).

Although the cell occupancy problem is well defined, there is no obvious connection between the modality of the density $f^n$ and the set of inequalities (5.2). We therefore again adopt the strategy of producing test densities derived from the taut string and gradually increase the modality until the inequalities (5.2) hold. The knowledge of which inequalities fail to hold provides further information which we are able to exploit as described in the next section.

5.2. *Local squeezing.* Local squeezing is described in Davies and Kovac (2001). We apply it to the density problem as follows. Suppose that one of the inequalities of (5.2) fails. We suppose that
$$w_{j,k}^n = |\{l : k2^{-j} < F^n(x_l) \leq (k+1)2^{-j}\}| \geq \upsilon_{j,k}^n(\alpha).$$
Clearly, there exists an interval $[x_{(l1)}, x_{(l2)}]$ such that $k2^{-j} < F^n(x_l) \leq (k+1)2^{-j}$ for all points $x_l$ in $[x_{(l1)}, x_{(l2)}]$. We now squeeze the tube locally on this interval and do this for all intervals where the upper inequality fails. For coefficients $w_{j,k}$ we proceed similarly but use slightly larger intervals such that $k2^{-j} < F^n(x_l) \leq (k+1)2^{-j}$ for all points $x_l$ in $(x_{(l1)}, x_{(l2)})$. The general procedure for doing this is as follows. First, a suitable initial global tube radius $\gamma_0$ is chosen using the Kolmogorov or generalized Kuiper metrics and the taut string is calculated. If all the inequalities (5.2) hold, the procedure terminates. If not, we reduce the radius by a factor $\rho$, $0 < \rho < 1$, on all intervals where an inequality fails. Typical choices for $\rho$ are 0.9 or 0.95. The taut string through the modified tube is calculated and, using this new test distribution, it is checked whether the inequalities (5.2) hold. If so, the procedure terminates. Otherwise the tube radius is again decreased by the factor $\rho$ on all intervals where an inequality fails. This is continued until all the inequalities are satisfied.

It is not easy to analyze the behavior of the local squeezing procedure. In the case of nonparametric regression Davies and Kovac (2001) give a heuristic argument indicating that the procedure improves the behavior at local extremes. A similar argument can be given for densities, but as it remains heuristic we omit it.

The ability of the local squeezing method to detect low power peaks [see Davies and Kovac (2001)] is shown by the following example. The data consist of a sample of size 1000 drawn from the four-normal distribution $N4$ of Section 3.1. The density is shown in the upper left corner of Figure 6. It exhibits a main peak, a moderate peak on the right and in the center two low power but very concentrated and very close peaks.

The upper right panel shows a kernel estimate which was calculated using a Gaussian kernel. The mode on the right-hand side was detected, but is considerably broader than the normal component of the original density function. The main component is well captured but there are three superfluous peaks. Finally, the two sharp peaks in the center of the data result in one flat local maximum. The lower left panel shows the result with the taut string method and two global tube radii. The solid line is derived from the $d_{\mathrm{KU}}^1$ metric. There are no spurious local extremes but the small central peaks are not detected. The dashed line shows that further global squeezing would only lead to additional spurious modes on the left before the central peaks are detected. Finally, the lower right panel shows the result of local squeezing. The number and locations of the local extrema are estimated correctly and the difference with respect to the original density function is very small.

Table 5 shows the performance of the local squeezing procedure for the distribtions $S$, $N1$, $N2$, $N4$, $N5$, $N10\_5$ and $N10\_10$ for samples of sizes 250 and 500. The procedure was calibrated as for the Kuiper metrics but, due to the discrete nature of the cell counts, it was not possible to adjust the parameters so that in 50% of the cases the modality for uniform samples was one. The choice lay between 48% and 55% and we took the latter. The results show a much enhanced performance for the distribution $N4$, but the results for the other distributions are worse than for the Kuiper metrics. This suggests a compromise procedure.

TABLE 5
*Results for the local squeezing procedure. The numbers give the percentage of simulations in which the correct modality was obtained. The numbers in parentheses give the mean absolute deviation from the correct modality. The results are based on* 5000 *simulations with sample sizes of* 250, 500 *and* 1000

| Dist. | $S$ | $N1$ | $N2$ | $N4$ | $N5$ | $N10\_5$ | $N10\_10$ |
|---|---|---|---|---|---|---|---|
| $n = 250$ | 91 (0.1) | 83 (0.2) | 42 (0.6) | 1 (1.6) | 4 (2.2) | 2 (2.9) | 99 (0) |
| $n = 500$ | 89 (0.1) | 80 (0.2) | 45 (0.6) | 22 (0.9) | 17 (1.5) | 36 (0.9) | 100 (0) |
| $n = 1000$ | 88 (0.1) | 79 (0.3) | 54 (0.5) | 75 (0.3) | 43 (0.8) | 91 (0.1) | 100 (0) |

5.3. *Compromise default procedures.* Statistical procedures make no assumptions about the data [Tukey (1993a)] and consequently are required to be compromises [see Tukey's example of the milk bottle in Tukey (1993b)]. Given a Kuiper metric $d_{\mathrm{KU}}^\kappa$, we calibrate the procedure based upon it so that in 60% of the cases the approximation to uniform samples is unimodal. Local squeezing is then applied so that the final approximation is unimodal in 50% of the cases. Again due to the discrete nature of the cell counts, 50% is not exactly attainable so we take the smallest percentage higher than 50. A second choice is to standardize the Kuiper procedure so that in 95% of the cases the approximation to uniform samples is unimodal. This is then reduced to 90% using local squeezing. We modify the local squeezing procedure as follows. Instead of using a multiresolution scheme we consider all intervals of length at most $\sqrt{n}$. This results in a procedure of $O(n^{1.5})$ but easily calculable for sample sizes of 50,000 and more. The reasoning behind this alteration is that we use local squeezing only to detect low power concentrated peaks. The others should be detected by the preceding Kuiper procedure. For reasons of space and comprehensibility we do not give an exact description of the local squeezing procedure but the source code is available from our web site. This leaves open the choice of $\kappa$ in $d_{\mathrm{KU}}^{\kappa}$. The software is available for all choices $\kappa = 1, 3, \ldots, 19$ with the default choice $\kappa = 19$. If data is to be analyzed in a routine manner, $\kappa$ can be chosen on the basis of experience or knowledge of the data involved.

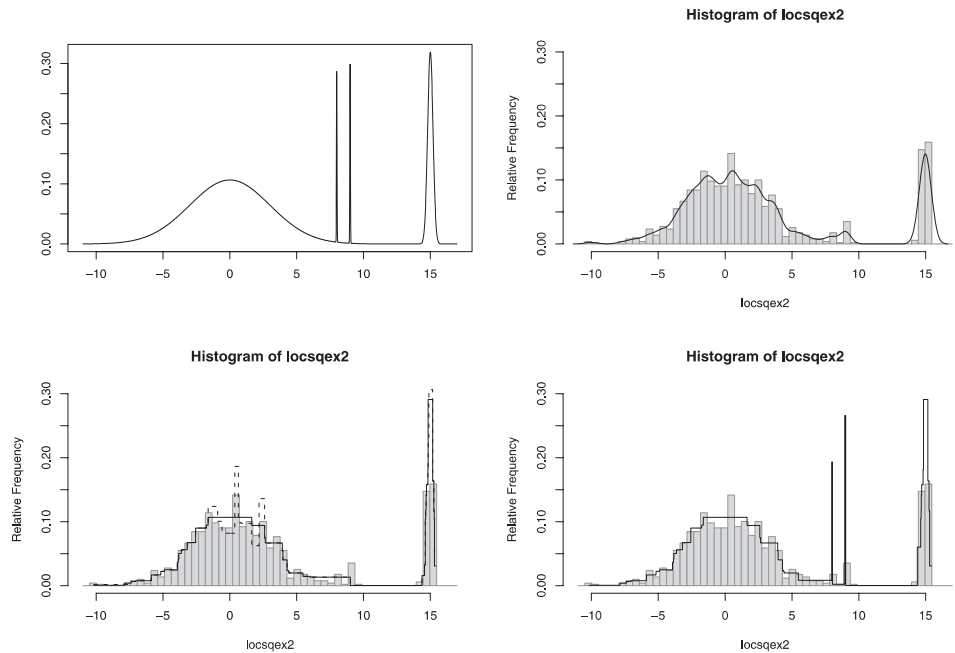

FIG. 6. *Local squeezing. The upper left panel shows the density of* $N4$. *A kernel estimate is shown in the upper right panel. The lower left panel illustrates global squeezing first with a solid line using the Kolmogorov bounds and then with a dashed line the taut string density with four modes. The local squeezing estimate is depicted in the lower right panel.*

TABLE 6
*Results for the compromise procedure based on $d_{\mathrm{KU}}^{19}$. The numbers give the percentage of simulations in which the correct modality was obtained. The numbers in parentheses give the mean absolute deviation from the correct modality. The results are based on* 1000 *simulations with sample sizes of* 250, 500 *and* 1000

| Dist. | $S$ | $N1$ | $N2$ | $N4$ | $N5$ | $N10\_5$ | $N10\_10$ |
|---|---|---|---|---|---|---|---|
| $n=250$ | 97 (0) | 93 (0.1) | 51 (0.5) | 2 (1.8) | 17 (1.6) | 40 (0.9) | 99 (0) |
| $n=500$ | 97 (0) | 94 (0.1) | 64 (0.4) | 19 (1.1) | 60 (0.5) | 95 (0) | 100 (0) |
| $n=1000$ | 99 (0) | 98 (0) | 86 (0.1) | 82 (0.2) | 99 (0) | 100 (0) | 100 (0) |

5.4. *Further simulations.* We now evaluate the two procedures COMPKU19_50 and COMPKU19_90 which are the compromise procedures described in the previous section using the Kuiper metric $d_{\mathrm{KU}}^{19}$ and calibrated at the uniform distribution to give the correct modality with probabilities 0.5 and 0.9, respectively. We compare them with two kernel-based methods. The first KERNCV uses likelihood cross-validation for the choice of the bandwidth, while the second KERNSJ uses the Sheather–Jones plug-in bandwidths. The comparisons are performed using the ten densities shown in Figure 7. They are taken from Marron and Wand (1992) and are the uniform distribution on $[0, 1]$, the Gaussian distribution and eight mixtures of normal distributions.

Each method was applied to 1000 samples of each of the densities and three different sample sizes (100, 500, 2000). For each estimate it was checked if the correct number of modes was found and if the positions of the modes corresponded to those of the densities. Table 7 shows how often the modes were determined correctly for the various densities and methods. Some comments are in order. First, if we use the procedure COMPKU5_50 which is tuned to three modes, then the performance for the trimodal density improves. For $n = 500$ three modal values are found in 20% of the cases, and for $n = 2000$ this rises to 37%. Second, all the densities are mixtures of a small number of Gaussian distributions with the exception of the uniform density for which the kernel methods based on a Gaussian kernel fail. The trimodal distribution is the one where the kernel methods perform clearly better than the taut string method. If, however, the central Gaussian distribution is replaced by a uniform distribution, then the kernel methods again fail. We refer to Hartigan (2000) for an explanation of this. It indicates that the comparison is weighted in favor of the kernel methods as both they and the densities are based on the Gaussian kernel. We note that the performance of the kernel methods seems to deteriorate with increasing sample size.

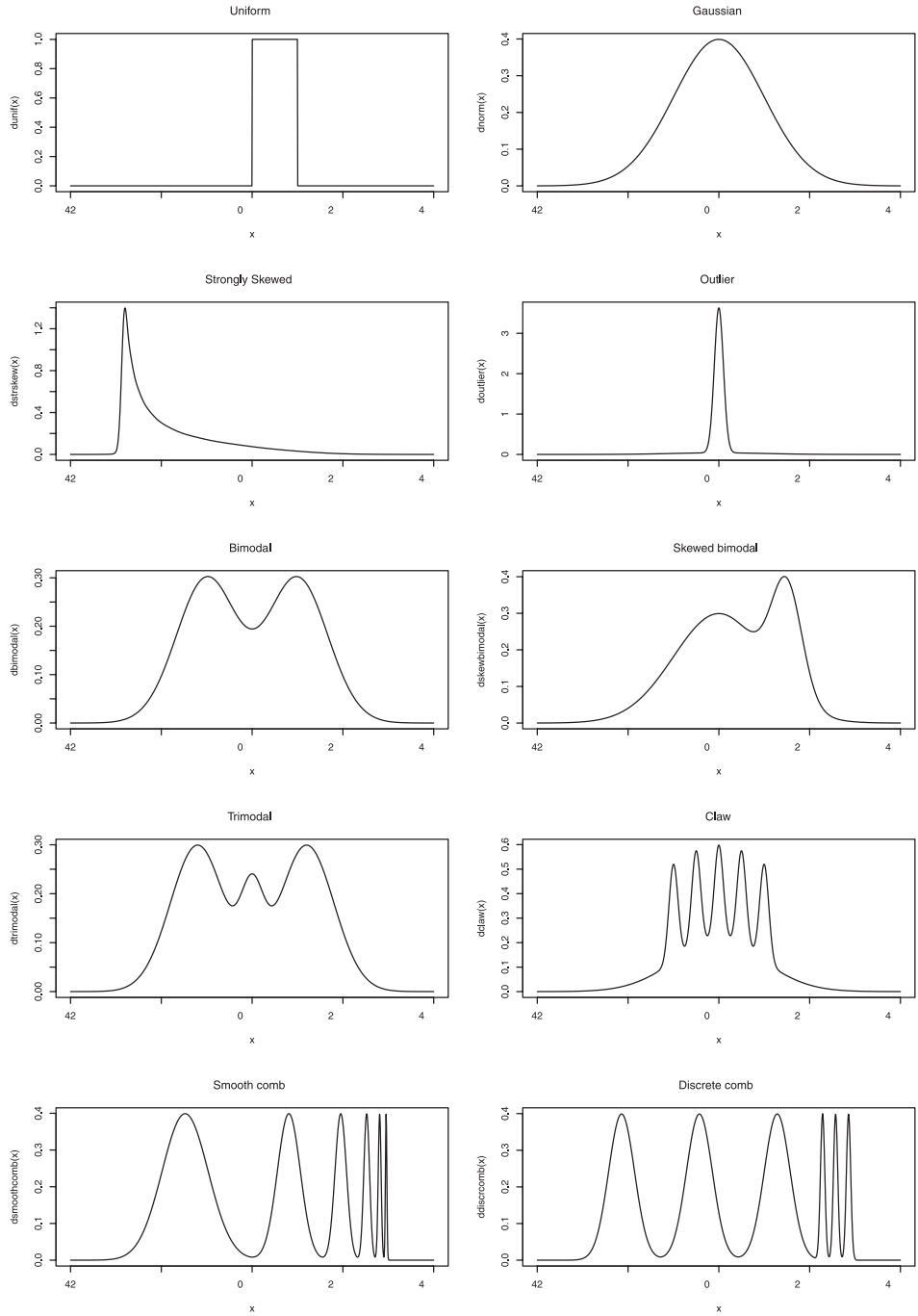

Fig. 7. *Ten densities that were used in a simulation study.*

## 6. Hidden periodicities, spectral densities and taut strings.

6.1. *Hidden periodicities.* The second problem we consider is that of detecting hidden periodicities in a data set $\mathbf{x}_n$. One method of formulating the problem is the following: calculate an appropriate spectral density function $f^n$ and identify the hidden periodicities in the data with the peaks of $f^n$ [Brillinger (1981), Priestley (1981) and Brockwell and Davis (1987), and the references given therein].

Existing methods by and large belong to one of two different categories of procedures. The first is nonparametric and uses some form of smoothing of the periodogram. This may take the form of kernel estimators or splines or wavelets or averages of periodograms obtained by splitting the data into blocks [see Brillinger (1981), Chapter 5, Neumann (1996) and the references given therein]. The second possibility is to model the data by an autoregressive process whose order is determined using some criterion such as AIC [Akaike (1977)], BIC [Akaike (1978)] or HQ [Hannan and Quinn (1979)]. The spectral density associated with the autoregressive process is then used to determine the hidden periodicities. None of these methods controls the number of peaks directly although the problem of hidden peaks is one of modality.

Before proceeding, we assume that the data have been normalized to have sample mean zero and variance 1. To ease the notation the transformed data will also be denoted by $\mathbf{x}_n$. In the context of time series $e_n$ will denote the empirical spectral density or the periodogram defined by

$$e_n(\omega) = \frac{1}{2\pi n} \left| \sum_{t=1}^{n} x_t \exp(i\omega\ t) \right|^2,$$

(6.1)
$$0 \leq \omega \leq 2\pi.$$

The corresponding empirical spectral distribution function $E_n$ is given by

(6.2) $$E_n(\omega) = \int_0^\omega e_n(\lambda)\, d\lambda.$$

TABLE 7
*Correctly detected modes in samples of various densities and for several automatic methods*

| Density | Size | KERNCV | KERNSJ | COMPKU19_50 | COMPKU19_90 |
|---|---|---|---|---|---|
| Uniform | 100 | 1 | 16 | 50 | 91 |
|  | 500 | 0 | 1 | 53 | 89 |
|  | 2000 | 0 | 0 | 53 | 91 |
| Gaussian | 100 | 77 | 79 | 85 | 98 |
|  | 500 | 79 | 78 | 95 | 99 |
|  | 2000 | 74 | 59 | 98 | 99 |
| Strongly skewed | 100 | 4 | 0 | 90 | 99 |
|  | 500 | 1 | 0 | 96 | 100 |
|  | 2000 | 0 | 0 | 99 | 99 |
| Outlier | 100 | 15 | 0 | 90 | 99 |
|  | 500 | 0 | 0 | 97 | 100 |
|  | 2000 | 0 | 0 | 98 | 100 |
| Bimodal | 100 | 71 | 81 | 45 | 14 |
|  | 500 | 75 | 84 | 68 | 33 |
|  | 2000 | 75 | 73 | 97 | 92 |
| Skewed bimodal | 100 | 32 | 46 | 34 | 9 |
|  | 500 | 45 | 37 | 35 | 13 |
|  | 2000 | 34 | 12 | 49 | 22 |
| Trimodal | 100 | 29 | 12 | 11 | 1 |
|  | 500 | 57 | 67 | 11 | 2 |
|  | 2000 | 81 | 82 | 20 | 6 |
| Claw | 100 | 1 | 0 | 4 | 0 |
|  | 500 | 2 | 2 | 63 | 34 |
|  | 2000 | 0 | 0 | 100 | 100 |
| Smooth comb | 100 | 18 | 0 | 1 | 0 |
|  | 500 | 5 | 0 | 5 | 1 |
|  | 2000 | 1 | 1 | 89 | 80 |
| Discrete comb | 100 | 12 | 0 | 1 | 0 |
|  | 500 | 2 | 0 | 31 | 13 |
|  | 2000 | 0 | 82 | 98 | 99 |

The candidate spectral densities we use are based on the taut strings $S_n$ through the Kolmogorov tubes centered at $E_n$. We assume that the taut string is constrained to go through $(0, L_n(0))$ and $(2\pi, E_n(2\pi)) = (2\pi, 1)$, where $L_n$ denotes the lower boundary.

One difference with respect to the i.i.d. model is the fact that the empirical spectral distribution function is defined for all $\omega$. In practice a grid must be chosen which, when analyzing the asymptotic behavior on test beds, becomes increasingly fine. We use the Fourier frequencies $\frac{2\pi j}{n}$, $j = 0, \ldots, n-1$, where the data have been augmented by zeros to produce a power of 2. Choosing a finer grid has had no effect on the data sets we have analyzed so far.

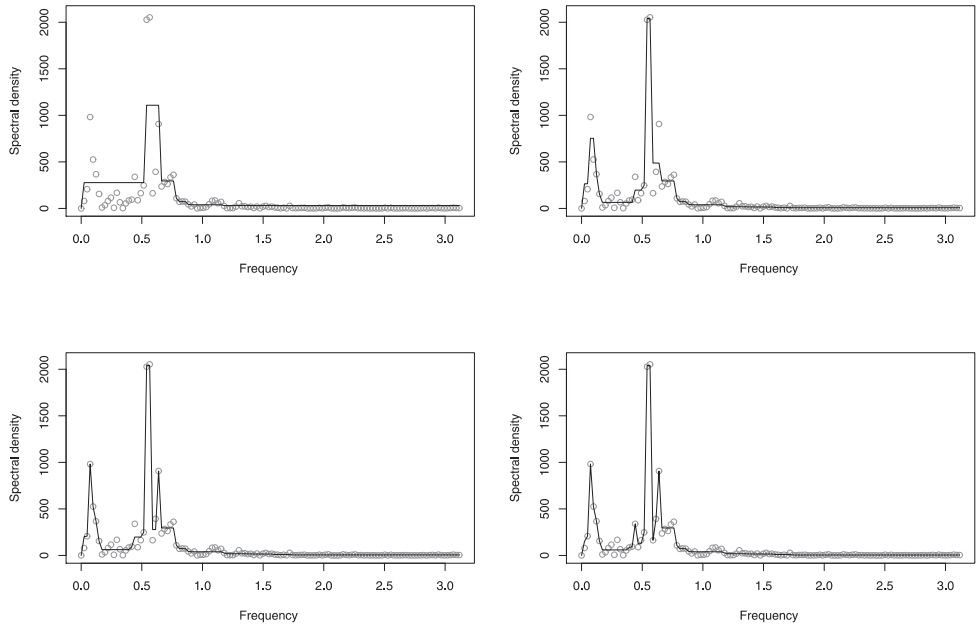

Fig. 8. *Sunspot data with number of peaks increasing from one to four.*

6.2. *Data analysis.* Just as in Section 3.3, it is possible to use the taut string as a data analytical tool. The radius of the Kolmogorov tube is gradually decreased and the resulting densities give information about the power and positions of the peaks. We give two examples. Figure 8 shows the first four peaks for the sunspot data [Anderson (1971)].

The second example is an artificial data set generated according to a scheme of Gardner (1988). Gardner does not explicitly specify the spectral density except that it has Gaussian shape with center frequency $2\pi\lambda$ with $\lambda = 0.35$. The density $f$ of (6.3) approximates the graph shown in Gardner's Figure 9.4(a) [Gardner (1988)]:

$$(6.3) \qquad f(\omega) = \tfrac{1}{3} e^{-300(\omega/2\pi - 0.35)^2}.$$

A realization of length 2048 was generated by filtering in the frequency domain. The following pure sine terms were added:

$$\sqrt{2}\sin(2\pi(0.2t - 106/360)),$$
$$\sqrt{2}\sin(2\pi(0.21t - 45.1/360)),$$
$$\sqrt{2}/10\sin(2\pi(0.1t - 32.6/360)).$$

A segment of length 256 starting at $t = 1023$ was taken as the simulated sample. It is shown in Figure 9.

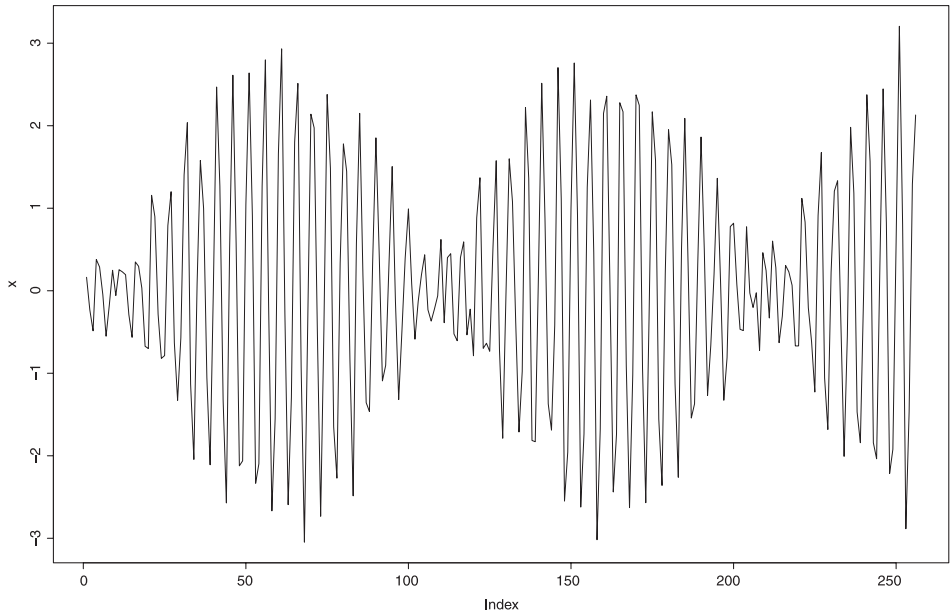

Fig. 9. *The Gardner data.*

A similar data set was analyzed by Gardner [(1988), Chapter 9.E, Experimental Study] in an experimental study of the performance of different spectral estimates. Figure 10 shows the first four peaks (in a log scale) for the data set of Figure 9. Finally, Figure 11 shows the four-peak density together with the periodogram.

6.3. *Two concepts of approximation.* The concepts of approximation used in the i.i.d. case had the advantage that the distributions involved were independent of the approximating model. This is no longer the case for stationary models. Furthermore, specifying the spectral distribution function $F$ does not specify the joint distribution of the stationary sequence. If, however, one is prepared to accept a Gaussian model, then the distribution $\mathbb{P}_F$ of the sequence is determined by $F$. In analogy with the i.i.d. case we have the following.

PROBLEM 6.1 (Kuiper spectral density problem). Determine the smallest integer $k_n$ for which there exists a spectral density $f^n$ with $k_n$ modes and whose distribution $F^n$ satisfies

(6.4) $$d_{\text{KU}}(E_n, F^n) \leq \text{qu}(n, \alpha, \mathbb{P}_{F^n}, d_{\text{KU}}),$$

where $\mathbb{P}_{F^n}$ denotes the distribution of the observations under the model.

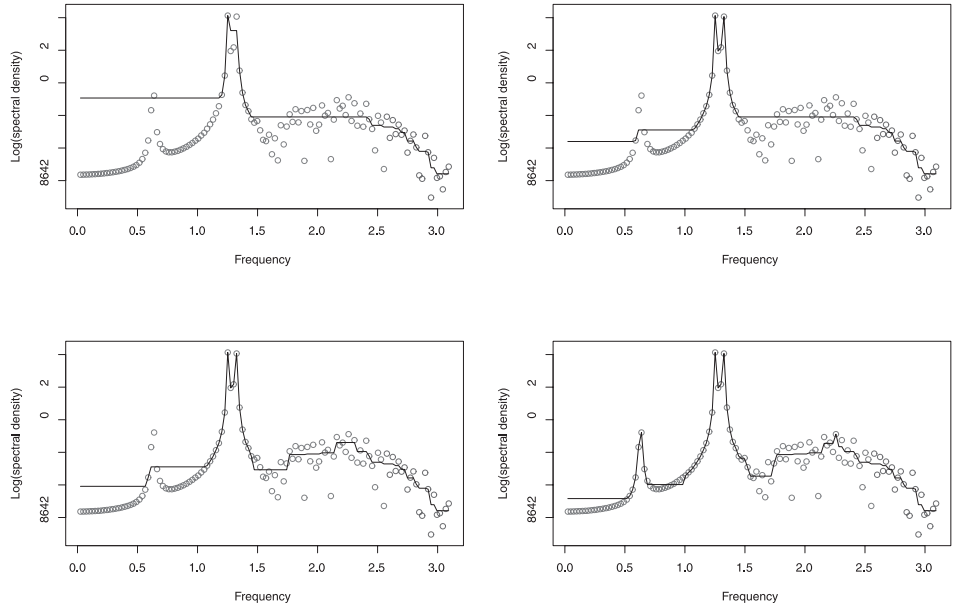

Fig. 10. *Gardner data with number of peaks increasing from one to four.*

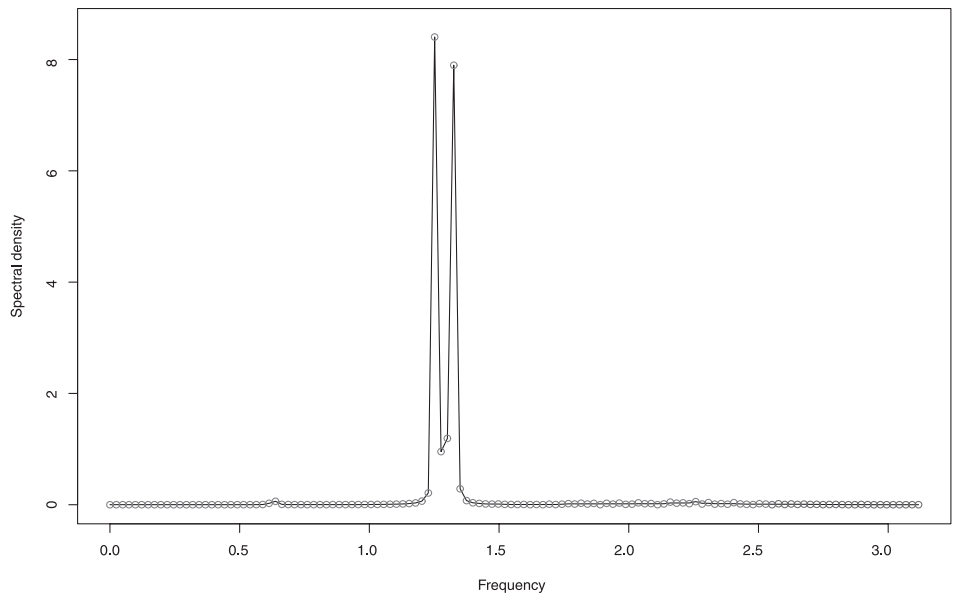

Fig. 11. *Gardner data with four peaks and the periodogram.*

There are two disadvantages with the procedure based on this concept of approximation. One is that the quantile in (6.4) depends on $F^n$. It would be possible to overcome this by using the taut string $S_n$ at each stage and then simulating the quantile $\mathrm{qu}(n, \alpha, \mathbb{P}_{S^n}, d_{\mathrm{KU}})$. This is clearly very time consuming. The second disadvantage is the following. Under appropriate conditions [Dahlhaus (1988)] we have the weak convergence result

$$\sqrt{n}(F_n - F) \Rightarrow Z,$$

where $F_n$ denotes the empirical spectral distribution function of the model with spectral distribution function $F$ and density $f$ and $Z$ denotes a continuous zero-mean Gaussian process defined by

(6.5) $$\mathbb{E}(Z(\lambda_1)Z(\lambda_2)) = \int_0^{\min(\lambda_1, \lambda_2)} f(\omega)^2 \, d\omega.$$

It follows from (6.5) that any large peaks will swamp smaller peaks which may be present and so prevent their detection. The one advantage of (6.4) is that it allows an asymptotic evaluation.

A more sensitive procedure is based on some kind of multiresolution analysis. Suppose for the moment that the sample size $n$ is a power of 2, $n = 2^m$. Given a spectral density function $f$, we define

(6.6) $$g_n(f, \omega) = \frac{e_n(\omega)}{f(\omega)},$$

and consider the multiresolution scheme

(6.7) $$w_{jk}(f) = \sum_{l=(j-1)2^k+1}^{j2^k} g_n(f, \omega_{l,n}),$$

$$j = 1, \ldots, 2^{m-k-1}, k = 0, \ldots, m-1,$$

where the $\omega_{l,n} = 2\pi l/n$ are the Fourier frequencies. The class of stationary processes with spectral density function $f$ is too large to provide a meaningful definition of approximation so we now restrict attention to Gaussian processes. Corresponding to level-dependent thresholds for wavelets, we specify lower and upper bounds $l_{k,n}$ and $u_{k,n}$, respectively, for the multiresolution coefficients (6.7). These now define the following.

PROBLEM 6.2 (Multiresolution spectral density problem). Determine the smallest integer $k_n$ for which there exists a spectral density $f^n$ with $k_n$ modes such that

(6.8) $\quad l_{k,n} \leq w_{jk}(f^n) \leq u_{k,n}, \qquad j = 1, \ldots, 2^{m-k-1}, k = 0, \ldots, m-1.$

The default bounds we use are $l_{k,n} = \text{qu}(\alpha_{1n}, 2^k)$ and $u_{k,n} = \text{qu}(\alpha_{2n}, 2^k)$, where $\text{qu}(\beta, \nu)$ denotes the $\beta$-quantile of the gamma distribution with $\nu$ degrees of freedom, $\alpha_{1n} = (1-\alpha)/2n$ and $\alpha_{2n} = 1 - \alpha_{1n}$ with $\alpha = 0.9$. The bounds are based on the Gaussian model and the asymptotic results for such processes as given, for example, by Theorem 5.2.6 of Brillinger (1981). If the asymptotic results hold precisely for finite $n$, then the bounds are chosen such that, for a stationary Gaussian process with spectral density function $f$, the inequalities (6.8) hold with probability at least 0.9 for $f^n = f$. As the individual $g_n(f, \omega)$ of (6.6) for $\omega = \frac{2\pi j}{n}$ are asymptotically independent, the bounds will be approximately of the correct order, again for Gaussian processes with a spectral density function. The usefulness of the bounds for real data sets is an empirical matter. In particular, they will be too slack if the spectral distribution function contains point masses.

This is the case for the Gardner data given above and may be seen in Figure 11. The absolutely continuous part of the spectrum shows a degree of noise, whereas the remainder of the spectrum is noise free. The default bounds we propose will detect the first peak but they are not sufficiently tight to split the two main peaks. On the other hand, if the bounds are sufficiently tight to separate the two peaks, then superfluous peaks will be produced in the absolutely continuous part of the spectrum. There would seem to be no easy solution which will work equally well for continuous as well as for discrete spectra.

We have no algorithm to solve the problem as it stands so again we use the local squeezing variant of the taut string method. The string is squeezed locally on the intervals where (6.8) fails and this is continued until all the inequalities are satisfied. When doing this, however, care must be taken regarding the order in which the inequalities are treated. From the form of $g_n(f, \omega)$ in (6.6) it is clear that a particular $g_n(f, \omega)$ can be very large and influence all intervals containing this particular frequency and this although the corresponding $e_n(\omega)$ is very small. Squeezing locally over all intervals affected by this frequency will often produce many superfluous peaks.

To avoid this we consider the intervals in order of size commencing with intervals of size 1. When all the inequalities are satisfied we then move on to intervals of size 2 and continue in this manner until all the inequalities are satisfied. This is the default version of the algorithm. If global squeezing is used, then the peaks will be introduced according to their power and may be introduced on intervals where the inequalities (6.8) are satisfied. This is the case for the Gardner data. If the default version with local squeezing is used, the main peak is not split. If, however, global squeezing is used, then it is split.

A practical problem which occasionally occurs is that the local squeezing version may find peaks of very small power which have no practical relevance. They may be removed by increasing the baseline of the empirical spectral density by a small amount. The software does this by first adding a small proportion of the total power, or the mean empirical spectral density, to the empirical spectral density and then proceeding as before.

6.4. *Asymptotics on test beds.* We indicate briefly the results of an asymptotic analysis using the Kuiper concept of approximation. The test bed we consider is that of a stationary process $X_n(F)$, $1 \leq n < \infty$, with a spectral distribution function $F$ and spectral density function $f$ as follows.

TEST BED 6.1.

(i) $F$ has exactly $k$ local extreme values on the interval $(0, \pi)$.
(ii) $F$ satisfies

$$\inf_{G \in \mathcal{F}(k-1)} \sup_{\omega \in [0,\pi]} |F(\omega) - G(\omega)| > 0,$$

where $\mathcal{F}(k-1)$ denotes the set of distributions with at most $k-1$ local extreme values.

To investigate the behavior of the taut string on the Test bed 6.1, we consider a tube of width $2C/\sqrt{n}$ and denote the taut string through this tube by $S_n(C)$ with derivative $s_n(C)$ and modality $k_n^C$. The intervals on which $s_n(C)$ takes on its local extreme values will be denoted by $I_i^e(n, C)$, $i = 1, \ldots, k_n^C$, with midpoints $\omega_i^e(n, C)$. The first theorem shows that on Test bed 6.1 the number and locations of the local extreme values are determined in a consistent manner.

THEOREM 6.1. *Consider the Test bed* 6.1. *Then for all* $\delta > 0$,

$$\lim_{C \to \infty} \liminf_{n \to \infty} \mathbb{P}\left(\{k_n^C = k\} \cap \left\{\max_{1 \leq i \leq k} |I_i^e(n, C)| \leq \delta\right\} \cap \left\{\max_{1 \leq i \leq k} |t_i^e(n, C) - t_i^e| \leq \delta\right\}\right)$$
$$= 1.$$

To obtain rates of convergence on appropriate test beds we must impose further conditions.

TEST BED 6.2.

(i) All spectral densities $f^j$ of order $j$ exist and $\sup_\omega |f^j(\omega)| \leq B^j$ for some constant $B$.
(ii) The spectral density function $f = f^2$ has a continuous second derivative $f^{(2)}$.
(iii) $f$ has exactly $k$ local extreme values, $0 < \omega_1, \ldots, \omega_k < 2\pi$, and $f^{(1)}(\omega) \neq 0$ for $\omega \in [0, 2\pi] \setminus \{\omega_1, \ldots, \omega_k\}$.

(iv) $f^{(2)}(\omega_j) \neq 0$, $j = 1, \ldots, k$.
(v) The fourth-order spectral density is continuous.

The above conditions correspond to (i) of Assumption 2.1 of Dahlhaus (1988).

Rates of convergence require a modulus of continuity for the process $Z_n = \sqrt{n}(F_n - F)$, where $F_n$ denotes the empirical spectral distribution function of the sample $(X_1(F), \ldots, X_n(F))$. Under the conditions of Theorem 2.4 of Dahlhaus (1988), it follows that

$$(6.9) \quad \sup_{0 \leq \omega_1 < \omega_2 \leq 2\pi, \omega_2 - \omega_1 < \delta} |Z_n(\omega_2) - Z_n(\omega_1)| \leq C\sqrt{\omega_2 - \omega_1} |\log(\omega_2 - \omega_1)|$$

with probability tending to 1 as $\delta$ tends to zero. From this it can be shown that the rate of uniform convergence away from the local extremes is $O(((\log n)^2/n)^{1/3})$. This differs from the rate of convergence for the test beds considered in Davies and Kovac (2001) by an extra $\log n$ term. This is explained by the different modulus of continuity. On the test beds of Davies and Kovac (2001) it is $\sqrt{\delta |\log \delta|}$, whereas above it is $\sqrt{\delta} |\log \delta|$.

6.5. *Examples.* The default version we use is the procedure deriving from the multiresolution problem with $\alpha = 1 - 0.1/n$ and a squeezing factor of 0.9. For the sunspot data the result is the one-peak density shown in the top left panel of Figure 8. For the Gardner data the result is the three-peak density derived from the four-peak density shown in the bottom right panel of Figure 10 but with the major peak not split (see above). Finally, we consider data generated according to a scheme of Neumann (1996), which is as follows:

$$X_n = Y_n + c_0 Z_n, \tag{6.10}$$

where

$$Y_n + a_1 Y_{n-1} + a_2 Y_{n-2} = b_0 \varepsilon_n + b_1 \varepsilon_{n-1} + b_2 \varepsilon_{n-2}$$

and $\{\varepsilon_n\}, \{Z_n\}$ are independent Gaussian white noise processes with variance 1. Neumann chose the coefficient values as follows: $a_1 = 0.2$, $a_2 = 0.9$, $b_0 = 1$, $b_1 = 0$, $b_2 = 1$ and $c_0 = 0.5$. A sample of size 1024 was generated according to this scheme. Figure 12 shows the logarithm of the spectral density of the sequence $\{X_n\}$ together with the logarithm obtained from the default version of the taut string method. The two peaks are correctly identified. The wavelet method used by Neumann results in six peaks [Neumann (1996), Figure 2(b)] for the data set he considered.

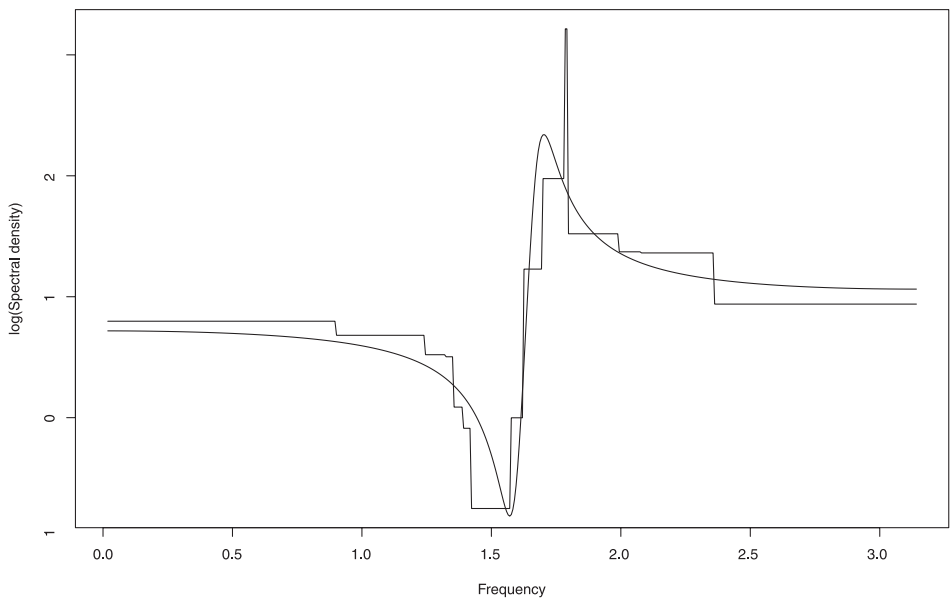

Fig. 12. *Log spectral densities of a sample of size* 1024 *generated by the scheme* (6.10).

## 7. Proofs.

7.1. *Proof of Theorem* 4.1. Using the Glivenko–Cantelli theorem, the property of the taut string of minimizing the modality in $T(F_n, \frac{C}{\sqrt{n}})$ and Proposition 12.3.3 of Dudley (1989) we see that

$$\min(\mathbb{P}(k_n^C \leq k), \mathbb{P}(k_n^C \geq k)) \geq \mathbb{P}\left(F \in T\left(F_n, \frac{C}{\sqrt{n}}\right)\right)$$

$$\geq 1 - \exp(-2C^2)$$

and conclude that

$$\lim_{C \to \infty} \lim_{n \to \infty} \mathbb{P}(k_n^C = k) = 1.$$

The other claims are proved similarly.

7.2. *Proof of Theorem* 4.2.

PROOF OF (a). Since the empirical process $E_n = \sqrt{n}(F_n - F)$ is tight, we conclude [Billingsley (1968), page 106] that

$$\lim_{C \to \infty} \lim_{n \to \infty} \mathbb{P}\left(\sup_{s \leq t \leq s + 2\tau_n} |E_n(s) - E_n(t)| \leq \frac{1}{C}\right) = 1,$$

where $\tau_n = \max(t_j^e - t_j^l)$, with $t_j^e$ denoting the point where $f$ takes its $j$th local extreme value and $t_j^l$ denoting the left endpoint of the $j$th local extreme interval of $f_n^C$.

From Theorem 4.1 we deduce that, for $C$ and $n$ sufficiently large, $f_n^C$ has the correct modality and

$$\sup_{s \le t \le 2\tau_n} |E_n(s) - E_n(t)| \le \frac{1}{C} \tag{7.1}$$

with arbitrarily high probability.

Suppose $F_n^C$ is initially convex and $t_1^l < t_1^e$. Then $F_n^C$ is the largest convex minorant of $F_n + C/\sqrt{n}$ [Barlow, Bartholomew, Bremner and Brunk (1972)] until it reaches the left endpoint $t_1^l(n,C)$ of $I_1^e(n,C) = [t_1^l(n,C), t_1^r(n,C)]$.

For some constant $\delta > 0$, for each $C$ and sufficiently large $n$,

$$t_1^r - t_1^e = \arg\max_{0 \le h \le \delta} H(h),$$

where

$$H(h) = \frac{F_n(t_1^l + h) - F_n(t_1^l) - 2C/\sqrt{n}}{h}. \tag{7.2}$$

As $F$ is convex on $[t_1^l, t_1^e]$, it can be shown using Taylor expansions that

$$G(h) = \frac{F(t_1^e + h) - F(t_1^e)}{h} \tag{7.3}$$

defines a strictly increasing function on $[0, \frac{4}{3}\mu]$, where $\mu = t_1^e - t_1^l$. Furthermore, for all $\tau < \mu$,

$$H\left(\frac{4}{3}\mu\right) - H(\tau) \ge G\left(\frac{4}{3}\mu\right) - G(\tau) + \frac{2C}{\sqrt{n\tau}} - \frac{2C}{\sqrt{n}4/3\mu} - \frac{2}{C\sqrt{n\tau}}$$
$$> 0.$$

This shows that $H$ cannot attain its maximum on $[0, \mu]$ and consequently $t_1^r > t_1^e$. Similar arguments hold for the other extrema. □

PROOF OF (b). We suppose that $S_n$ has a local maximum on $I_1^e(n,C) = [t_1^l(n,C), t_1^r(n,C)]$, that $t_1^e \in I_1^e$ and that (7.1) is satisfied. Define $G$ by

$$G(h) = \frac{F(t_1^l + h) - F(t_1^l) - 2C/\sqrt{n}}{h},$$

and consider $h_0 = \arg\max_{0 \le h \le \delta} G(h)$. Then $G'(h_0) = 0$ implies

$$f(t_1^l + h_0)h_0 = F(t_1^l + h_0) - F(t_1^l) - \frac{2C}{\sqrt{n}}.$$

Using Taylor expansions in $t_1^e$ and the fact that $f'(t_1^e) = 0$, we obtain

$$h_0^3 \ge -\frac{6C}{\sqrt{n}f''(t_1^e)} + o(h_0^3).$$

In the other direction we consider

(7.4) $$h_1 = \arg\max_{0 \leq h \leq \delta} \frac{F(t_1^e + h) - F(t_1^e) - 2C/\sqrt{n}}{h}$$

and

$$h_2 = \arg\min_{0 \leq h \leq \delta} \frac{F(t_1^e - h) - F(t_1^e) - 2C/\sqrt{n}}{h}.$$

It is not difficult to see that $h_0 \leq h_1 + h_2$. Setting the derivative of the right-hand side of (7.4) to zero and using a Taylor expansion in $t_1^e$ yields

$$h_1^3 = -\frac{6C}{\sqrt{n} f''(t_1^e)} + o(h_1^3).$$

The same argument holds for $h_2$ as well and both together show that

$$h_0^3 \leq -\frac{12C}{\sqrt{n} f''(t_1^e)} + o(h_0^3).$$

Define $H$ as in (7.2) and consider

$$\tilde{h}_0 = \arg\max G(h) - \frac{2}{\sqrt{Cnh}}.$$

The considerations above show that

$$\left( -\frac{6(C + 1/\sqrt{C})}{\sqrt{n} f''(t_1^e)} \right)^{1/3} \leq \tilde{h}_0(1 + o(1)) \leq \left( -\frac{12(C + 1/\sqrt{C})}{\sqrt{n} f''(t_1^e)} \right)^{1/3}.$$

Furthermore, considerations as in (a) show that $G(x) - \frac{2}{\sqrt{Cn}}$ defines a strictly decreasing function. Therefore, for all $h > (1 + \frac{1}{\sqrt{C}})\tilde{h}_0$,

$$H(\tilde{h}_0) - H(h) \geq G(\tilde{h}_0) - G(h) - \frac{2}{C\sqrt{nh}} > 0.$$

Consequently, $H$ cannot attain its maximum in $h > \tilde{h}_0(1 + \frac{1}{\sqrt{C}})$ and hence

$$\arg\max_{0 < h < \delta} H(h) < \left(1 + \frac{1}{\sqrt{C}}\right) \left( -\frac{12(C + 1/\sqrt{C})}{\sqrt{n} f''(t_1^e)} \right)^{1/3}.$$

Similarly, it can be shown that

$$\arg\max_{0 < h < \delta} H(h) < \left(1 - \frac{1}{1 + \sqrt{C}}\right) \left( -\frac{6(C - 1/\sqrt{C})}{\sqrt{n} f''(t_1^e)} \right)^{1/3}. \qquad \Box$$

PROOF OF (c). The proof relies on the modulus of continuity of the empirical process $E_n$.

LEMMA 7.1. *Let $Y(n,C)$ denote random variables such that, for all $\varepsilon > 0$,*

$$\lim_{C\to\infty}\lim_{n\to\infty} \mathbb{P}(|Y(n,C)| < \varepsilon) = 1.$$

*Consider $\alpha_n = n^{-\gamma}$ for some $\gamma < 1$ and*

$$\beta_n^C = \max\left\{\frac{1}{\log(n)}, Y(n,C)\right\}.$$

*Then for all $B > 2$ we have*

$$\lim_{C\to\infty}\lim_{n\to\infty} \mathbb{P}\left(\max_{\alpha_n < |s-t| < \beta_n^C} \frac{|E_n(s) - E_n(t)|}{\sqrt{|t-s|\log(1/|t-s|)}} > B\right) = 0.$$

PROOF. Define random integer-valued variables $K_n$ by

$$K_n = \left\lfloor \log_2\left(\frac{\beta_n^C}{\alpha_n}\right)\right\rfloor.$$

Using a result of Mason, Shorack and Wellner (1983), we conclude that provided $\beta_n^C < \frac{1}{2}$,

$$\mathbb{P}\left(\max_{\alpha_n < |s-t| < \beta_n^C} \frac{|E_n(s) - E_n(t)|}{\sqrt{|t-s|\log(1/|t-s|)}} > B\right)$$

$$\leq \sum_{k=0}^{\infty} \mathbb{P}\left(|E_n(s) - E_n(t)| > B\sqrt{|t-s|\log\left(\frac{1}{|t-s|}\right)}\right. \text{ for some } s,t \text{ with}$$

$$\left. 2^k\alpha_n < |s-t| < 2^{k+1}\alpha_n | k \leq K_n\right)$$

$$\leq \sum_{k=0}^{\infty} \frac{20}{a_k(\beta_n^C)^3} \exp\left(-(1-\beta_n^C)^4 \frac{\lambda_k^2}{a}\psi\left(\frac{\lambda_k}{\sqrt{na_k}}\right)\right),$$

where we denote $2^{k+1}\alpha_n$ by $a_k$,

$$\lambda_k = B\sqrt{\frac{\log(1/\alpha_n)}{2}}$$

and

$$\psi(x) = 2\frac{(1+x)(\log(1+x)-1)+1}{x^2}.$$

It is easily verified that $\psi(\frac{\lambda_k}{\sqrt{na_k}}) \to 1$. Thus

(7.5) $$\lim_{C,n\to\infty} \mathbb{P}\left((1-\beta_n^C)^4\psi\left(\frac{\lambda_k}{\sqrt{na_k}}\right) > \frac{2}{B}\right) = 1.$$

Putting this together, we deduce that

$$\mathbb{P}\bigg(|E_n(s) - E_n(t)| > B\sqrt{|t-s|\log\Big(\frac{1}{|t-s|}\Big)} \text{ for some } s, t \text{ with}$$

$$\alpha_n < |s - t| < \beta_n^C\bigg)$$

$$< \frac{20\log(n)^3}{n^{\gamma(B/2-1)}}.$$

This completes the proof of the lemma. $\square$

We proceed now with the proof of (c). Since $f$ is twice continuously differentiable, there is some constant $D > 0$ such that

$$|F(x+h) - F(x) - hf(x) - \tfrac{1}{2}h^2 f'(x)| \le Dh^3$$

for all $x$ and $h$.

Let $B$ be an arbitrary constant greater than 2 and

$$d(n,C) = \min\bigg\{|f'(x)| \,|\, x \in [0,1] \setminus \bigcup_i I_i^e(n,C)\bigg\}.$$

Define a random sequence $h(n,C)$ by

$$h(n,C) = \frac{(8B)^{2/3} \log(d(n,C)^2 n)^{1/3}}{(3n)^{1/3} d(n,C)^{2/3}}.$$

We consider the situation where

(i) $f_n^C$ attains the correct modality;
(ii) $t_i^e \in I_i^e(n,C)$ for all $i$;
(iii) the empirical process satisfies

$$\sup_{|s-t| < Y(n,C)} |E_n(t) - E_n(s)| < B\sqrt{|s-t|\log\Big(\frac{1}{|s-t|}\Big)},$$

where $Y(n,C)$ is defined by

$$Y(n,C) = \max\{x_{j+1} - x_j \,|\, x_j, x_{j+1} \text{ knots}, [x_j, x_{j+1}] \ne I_i^e(n,C) \text{ for all } i\};$$

(iv) for all $x \in [0,1] \setminus \bigcup_i I_i^e(n,C)$,

$$h_n \le \frac{f'(x)}{32D}$$

holds;
(v) for each extreme interval $I_i^e(n,C)$, the distances of each endpoint to $t_1^e$ are both smaller than $4h_n$.

The preceding lemmas and parts of this theorem show that the probability that all these assumptions are satisfied simultaneously converges to 1 as $n$ and $C$ tend to $\infty$. For example, (7.2) follows from (b) which provides a constant $A > 0$ such that $|f'(x)| \geq An^{-1/6}$.

Consider now an arbitrary point $t_1 \in [0,1] \setminus \bigcup_i I_i^e(n,C)$, where $f'(t_1) > 0$. Then

$$\frac{F_n(t_1 + h_n) - F_n(t_1)}{h_n} \leq f(t_1) + \frac{1}{2} h_n f'(t_1) + D h_n^2 + \frac{B\sqrt{\log(1/h_n)}}{\sqrt{nh_n}}.$$

Plugging in the expression for $h_n$ and using the assumptions made above, we see that

$$\frac{F_n(t_1 + h_n) - F_n(t_1)}{h_n} \leq f(t_1) + \frac{1}{2} h_n f'(t_1)\left(1 + \frac{1}{4} + \frac{1}{4}\right).$$

Similarly, we conclude that, for all $h \in [4h_n, t_j^e]$,

$$\frac{F_n(t_1 + h) - F_n(t_1)}{h} \geq f(t_1) + \frac{1}{2} h f'(t_1)\left(1 - \frac{1}{4} - \frac{1}{4}\right),$$

where $t_j^e$ is the smallest local extreme value greater than $t_1$.

Suppose that there are knots $x_j$ and $x_{j+1}$ that do not embrace a local extreme interval such that $h_0 = x_{j+1} - x_j > 4h_n$ and such that $f$ is increasing on $[x_j, x_{j+1}]$. The width $\tilde{h}$ is the local argmin

$$\tilde{h} = \arg\min_{0 < h < \delta} \frac{F_n(x_1 + h) - F_n(x_1)}{h}.$$

On the other hand, the considerations above show that

$$\frac{F_n(x_1 + h_n) - F_n(x_1)}{h_n} < \frac{F(x_1 + h) - F_n(x_1)}{h}.$$

Therefore, the distance between two knots that do not embrace an extreme interval is bounded by $4h_n$. □

PROOF OF (d). We assume that all the assumptions made in the proof of (c) are again satisfied and that each two extreme intervals $I_i^e$ and $I_{i+1}^e$ are separated by at least two additional knots $x_j$ and $x_{j+1}$:

$$\max I_i^e < x_j < x_{j+1} < \min I_{i+1}^e.$$

Define $h_n$ as in (7.4). Consider a knot $x_i$ which does not delimit a local extreme interval $I_i^e$. We take $f$ to be increasing in $x_i$. Then the proof of (c) shows that

$$f_n^C(x_i) \leq \frac{F_n(x_i + h_n) - F_n(x_i)}{h_n}$$
$$\leq f(x_i) + C_1 |f'(x_i)|^{1/3} \left(\frac{\log(n)}{n}\right)^{1/3}.$$

Similar arguments show that

$$f_n^C(x_i) \geq \frac{F_n(x_i) - F_n(x_i - h_n)}{h_n}$$

$$\geq f(x_i) - C_1|f'(x_i)|^{1/3}\left(\frac{\log(n)}{n}\right)^{1/3}.$$

Analogous inequalities can be derived in the case where $f$ is decreasing in $x_i$.

Suppose now that $t$ is an arbitrary point in

$$\left[A\left(\frac{\log(n)}{n}\right)^{1/3}, 1 - A\left(\frac{\log(n)}{n}\right)^{1/3}\right] \setminus \bigcup_{i=1}^{k} I_i^e(n, C).$$

Let $x_i$ be the nearest knot which does not delimit a local extreme interval. Then

(7.6)
$$|f(t) - f_n^C(t)|$$
$$\leq |f(t) - f(x_i)| + |f(x_i) - f_n^C(x_i)| + |f_n^C(x_i) - f_n^C(x)|.$$

The inequalities above show that the second term is bounded by

$$C_2|f'(x_i)|^{1/3}\left(\frac{\log(n)}{n}\right)^{1/3}.$$

The first term is bounded by

$$C_3|t - x_i|\,|f'(x_i)| \leq C_3|f'(x_i)|^{1/3}\left(\frac{\log(n)}{n}\right)^{1/3}.$$

This follows from (b).

Depending on the exact definition of $f_n^C(x)$ at knot points, the third term is either 0 or bounded by $2C_1|f'(x_i)|^{1/3}(\frac{\log(n)}{n})^{1/3}$.

This completes the proof of (d). □

PROOF OF (e). As in the other cases, we assume that $f_n^C$ attains the correct modality and that $t_i^e \in I_i^e(n, C)$ for each extreme point $t_i^e$. We also assume that, for each extreme interval $I_i^e$,

$$\left(1 - \frac{1}{1+\sqrt{C}}\right)\left(-\frac{6(C - 1/\sqrt{C})}{\sqrt{n}f''(t_1^e)}\right)^{1/3}$$
$$\leq |I_i^e(n, C)|$$
$$\leq \left(1 + \frac{1}{\sqrt{C}}\right)\left(-\frac{12(C + 1/\sqrt{C})}{\sqrt{n}f''(t_1^e)}\right)^{1/3}.$$

The regression function $f_n^C$ takes in $t_i^e$ the slope of the taut string in the extreme interval $I_i^e = [x_1, x_2]$. Taylor expansions in $t_i^e$ using $f'(t_i^e) = 0$ and an application of the modulus of continuity for the empirical process $E_n$ as formulated in Lemma 7.1 yield

$$|f_n^C(t_i^e) - f(t_i^e)| \leq D_1(1 + o(1))\frac{f''(t_i^e)^{1/3}}{n^{1/3}}.$$

The proof is now completed by extending the bound to arbitrary points in extreme intervals $I_i^e$. This is done in the usual way as in (7.6) using a Taylor expansion in $t_i^e$ and shows that

$$|f(t) - f(t_i^e)| \leq D_2 |I_i^e|^2 f''(t_i^e). \qquad \square$$

**Software.** The software is available from our home page at www.stat-math.uni-essen.de. A package for R is in preparation.

**Acknowledgments.** We thank two referees and an Associate Editor for their comments on the first version of the paper. The final version has profited from their remarks and suggestions.

**Note added in proof.** After acceptance of this paper for publication we found that a small change of the notion of adequacy for densities leads to a considerable improvement in the performance of the procedure. In particular,

- a calibration using the uniform density is not necessary;
- a constant density is fitted for almost all samples of the uniform distribution;
- the peaks of densities such as the claw density are detected more reliably.

We consider differences of Kuiper metrics $d_{\text{KU}}^\kappa$

$$(*) \qquad \rho_i(F, F_n) = d_{\text{KU}}^i(F, F_n) - d_{\text{KU}}^{i-1}(F, F_n), \qquad i = 1, \ldots, \kappa,$$

where $d_{\text{KU}}^0 \equiv 0$. The distribution of each difference $\rho_i(F, F_n)$ is independent of $F$. In our modified $\kappa$-Kuiper problem we require all differences to be smaller than some $\alpha$-quantile of $\rho_i$ with $\alpha$ close to 1. Our default is $\alpha = 0.999$.

PROBLEM 3.2' (Modified $\kappa$-Kuiper density problem). Determine the smallest integer $k_n$ for which there exists a density $f^n$ with $k_n$ modes and whose distribution $F^n$ satisfies

$$\rho_i(E_n, F^n) \leq \text{qu}(n, \alpha, \rho_i)$$

for all $i = 1, \ldots, \kappa$.

TABLE 8
*Results for the taut string procedure based on the modified $\kappa$-Kuiper criterion using $\kappa 19$. The numbers give the percentage of simulations in which the correct modality was obtained. The numbers in parentheses give the mean absolute deviation from the correct modality. The results are based on* 1000 *simulations with sample sizes of* 100, 250, 500 *and* 1000

| Dist. | $U$ | $S$ | $N1$ | $N2$ | $N4$ | $N5$ | $N10\_5$ | $N10\_10$ |
|-------|---------|----------|----------|---------|--------|----------|---------|---------|
| 100   | 99 (0.0)| 100 (0.0)| 98 (0.0) | 8 (1.2) | 0 (2.1)| 3 (3.8)  | 23 (3.8)| 63 (0.4)|
| 250   | 99 (0.0)| 100 (0.0)| 100 (0.0)| 18 (0.8)| 0 (2.0)| 23 (2.6) | 80 (0.2)| 91 (0.1)|
| 500   | 98 (0.1)| 100 (0.0)| 100 (0.0)| 53 (0.5)| 1 (2.0)| 76 (0.4) | 94 (0.1)| 95 (0.1)|
| 1000  | 98 (0.1)| 100 (0.0)| 100 (0.0)| 86 (0.1)| 3 (2.0)| 100 (0.0)| 98 (0.0)| 97 (0.0)|

As for the original $\kappa$-Kuiper problem quantiles may again be obtained by simulation. For large $n > 1000$ the distribution of $\sqrt{n}\rho_i(F_n, F)$ can be approximated by the corresponding quantile of a Brownian bridge using the weak convergence of the empirical process.

The taut string procedure can be initiated by using the global bandwidth $\varepsilon^0 = 0.5$ which corresponds to a constant approximating density $f^0$. If all inequalities $(*)$ are satisfied for $f^0$ we are finished. Otherwise assume that $i$ is the smallest index such that an inequality $(*)$ is not satisfied. Then we set $\varepsilon^1 = 0.5 \cdot \text{qu}(n, \alpha, \rho_i)$ which is the largest tube width such that the $i$-th difference $\rho_i$ of $\kappa$-Kuiper metrics is sufficiently small. After a few iterations all $\rho_i$ will satisfy $(*)$ and the final approximation is the taut string approximation with the maximal global bandwidth and hence minimal number of modes which is adequate for the data.

Table 8 shows that the proposed procedure returns a constant function for the uniform density in about 99 per cent of the cases independent of $n$. At the same time the 10 peaks of the $N10\_5$-density are found in 23 % of the cases by the new procedure for samples of size 100 and in 80 % for samples of size 250. The old procedure never found the correct number of peaks for samples of size 250. Only the performance for the bimodal $N2$-distribution has deteriorated. The small peaks of the $N4$ distribution are only detected occasionally even for large sample sizes.

Universität Duisburg-Essen
Campus Essen
Fachbereich Mathematik
45117 Essen
Germany
e-mail: laurie.davies@uni-essen.de
e-mail: Arne.Kovac@uni-essen.de